\documentclass[11pt]{article}

\usepackage[T1]{fontenc}
\usepackage[utf8]{inputenc}
\usepackage{newtxtext}
\usepackage[a4paper,margin=2.4cm]{geometry}
\usepackage[authoryear,round]{natbib}
\usepackage{amsmath,amssymb}
\usepackage{graphicx}
\usepackage{booktabs}
\usepackage{array}
\usepackage{xcolor}
\usepackage{tikz}
\usetikzlibrary{arrows.meta,positioning}
\usepackage[font=small,labelfont=bf,labelsep=period]{caption}
\usepackage{authblk}
\usepackage{url}
\usepackage[colorlinks=true,linkcolor=blue!50!black,citecolor=blue!50!black,urlcolor=blue!50!black]{hyperref}

\usepackage{etoolbox}
\AtBeginEnvironment{table}{\renewcommand\scriptsize{\fontsize{7}{8}\selectfont}}
\newcommand{\tabnotemark}[1]{\textsuperscript{#1}}
\newcommand{\tablenote}[1]{\par\vspace{3pt}{\footnotesize #1\par}}
\title{Agent-written heuristics versus commercial solvers on unsolved job-shop benchmark instances}

\author[1,*]{Andre Albrecht}
\author[1]{Enzo Morosini Frazzon}
\author[2]{Leandro C. Coelho}
\author[1,3]{Bruno Vieira}
\affil[1]{ProLogIS -- Production and Logistics Intelligent Systems, Department of Production Engineering, Federal University of Santa Catarina (UFSC), Florian\'opolis, SC, 88040-900, Brazil}
\affil[2]{CIRRELT and Canada Research Chair in Integrated Logistics, Universit\'e Laval, Qu\'ebec, QC, G1V 0A6, Canada}
\affil[3]{Department of Industrial Engineering, Federal Institute of Santa Catarina (IFSC), Ca\c{c}ador, SC, 89503-550, Brazil}
\affil[ ]{\texttt{andre.albrecht@posgrad.ufsc.br}, \texttt{enzo.frazzon@ufsc.br}, \texttt{leandro.coelho@fsa.ulaval.ca}, \texttt{bruno.vieira@ifsc.edu.br}}
\affil[*]{Corresponding author}

\hypersetup{pdftitle={Agent-written heuristics versus commercial solvers on unsolved job-shop benchmark instances},pdfauthor={Andre Albrecht, Enzo Morosini Frazzon, Leandro C. Coelho, Bruno Vieira}}

\date{\small Preprint submitted to \emph{International Transactions in Operational Research}, November 2026}

\begin{document}

\maketitle

\begin{abstract}
Commercial solvers are the standard tool for difficult scheduling problems. We study whether a frontier large language model, used as a coding agent with open-source tools only, can match them in quality and speed. We report every run, and independent checkers verify every schedule. On the 100{,}000-operation job shop \texttt{tai\_100\_1000\_1}, an agent-written solver reached 525{,}515, below OptalCP, Hexaly, and CP Optimizer. Unchanged, it improved all ten instances of the family and closed two. With equal budgets on the same workstation, it ended 0.0--1.8\% above the lower bound, against 2.6--5.8\% for CP Optimizer and 2.4--4.4\% for Hexaly. On the Behnke--Geiger flexible job shops, an agent-written tabu search improved 11 of 32 open instances. With two transport robots, an agent-written solver improved seven of ten MKT instances. New lower bounds show that four published MKT values cannot be attained on the distributed instance files. Agent-run improvement campaigns closed eight further instances.
\end{abstract}


\section{Introduction}\label{sec:intro}

Scheduling problems arise in manufacturing and logistics \citep{leusin2018jssp}, and many practitioners solve them by means of commercial solvers. Typical examples are Hexaly \citep{hexaly,benoist2011localsolver}, IBM ILOG CP Optimizer \citep{laborie2018}, Gurobi, and OptalCP \citep{optalcp}. These solvers perform very well on a wide range of problems and need little problem-specific work, but commercial use requires a paid license. Hexaly and OptalCP quote their commercial prices on request \citep{hexalypricing,optalcpeditions}, and IBM lists a developer subscription of CP Optimizer from US\$320 per user and month, with deployment licenses on request \citep{ibmcplexpricing}. These three are free for academic use. Open-source solvers such as OR-Tools CP-SAT \citep{ortools,perron2023cpsatlp} are free, but on large scheduling instances they are often weaker \citep{dacol2022,yuraszeck2026}.

At the same time, large language models (LLMs) can now write and run code in long autonomous sessions \citep{chen2021codex,jimenez2024swebench,yang2024sweagent}. Recent work uses LLMs to design heuristics \citep{romera2024funsearch,liu2024eoh,ye2024reevo,novikov2025alphaevolve} and to formulate optimization models \citep{ramamonjison2023nl4opt,ahmaditeshnizi2024optimus,huang2025orlm}. Scheduling is a natural test for such agents, because a candidate schedule is easy to check, whereas a good schedule is difficult to find.

Our starting question is whether a frontier LLM, used as a coding agent, can replace an off-the-shelf commercial solver. We aimed for results of similar quality in a similar time. The agent had no access to a commercial solver, and it could use only open-source tools and the code it wrote itself.

The results exceed this target. The agents analyzed the structure of each instance and, when the standard formulation did not yield good solutions, they changed it. They then wrote heuristics and metaheuristics tailored to that structure, and used open-source tools such as CP-SAT for subproblems and bounds. On the job shop and the flexible job shop, the results are better than the best published values of the commercial solvers, and better than CP Optimizer and Hexaly in our own runs with the same budget on the same workstation. For the transport problem, the results are better than all valid published values on seven of the ten MKT instances, including the largest ones, better than the CP Optimizer models of \citet{meng2025cp} on nine of them, and better than CP-SAT with four times the CPU time. We do not count four published values on these instances, because they are below a new lower bound on the distributed instance files, nor two published schedules that assign an operation to a machine on which it is not eligible. On the Chambers--Barnes transport instances, where processing times dominate travel times, the CP Optimizer approach is better than the agent-written solver.

The contributions of the paper can be summarized as follows:
\begin{itemize}
\item We propose a method to test LLM coding agents on open benchmark instances, with a clear split between autonomous runs and runs with human help (Section~\ref{sec:method}).
\item We provide new best-known schedules for three established problem families: the $100\times1000$ job shops of \citet{dacol2022}, the flexible job shops of \citet{behnke2012}, and the flexible job shops with transport robots of \citet{homayouni2021} (Section~\ref{sec:results}).
\item We derive new lower bounds for the flexible job shop with transport robots and validate them with two solvers. They are above eleven published values on the distributed instance files and certify eight published schedules as optimal, two of them (MFJST09 and \texttt{mt10xyzt}) without a previous proof (Sections~\ref{sec:fjsptlb} and~\ref{sec:res-fjspt}).
\item We compare the quality, speed, and cost of the agent-written solvers with those of commercial and open-source solvers on the same instances (Sections~\ref{sec:results} and~\ref{sec:disc}), and three ways to use an agent: as the solver itself, as the author of a reusable solver, or as the improver of such a solver in its own campaign (Section~\ref{sec:modes}).
\item We test the agent-written solvers on similar instances without any change. On the job-shop family, the solver has improved the best-known values of all ten instances and has solved two open instances to optimality (Section~\ref{sec:res-jssp}).

\end{itemize}

We report every agent run on the instances of this paper, without any selection based on the results. An independent checker verifies every schedule, and a public repository contains all instances, schedules, checkers, solvers, solver logs, and task texts (see the statement on data availability).

The remainder of the paper is organized as follows. In Section~\ref{sec:related} we review the literature on shop scheduling, on the flexible job shop with transport robots, and on LLMs for optimization. The three problems and new lower bounds for the problem with transport robots are formally described in Section~\ref{sec:problem}, and their selection in Section~\ref{sec:selection}. In Section~\ref{sec:method} we present our experimental method, followed by extensive computational experiments in Section~\ref{sec:results} and a discussion in Section~\ref{sec:disc}. Conclusions follow in Section~\ref{sec:conc}, and Appendix~\ref{app:newbest} lists all new best-known values.

\section{Literature review}\label{sec:related}

\subsection{Solvers and benchmarks for shop scheduling}

The job shop and the flexible job shop are classical NP-hard problems \citep{garey1976,brandimarte1993}, and methods are compared on public benchmark libraries \citep{taillard1993,beasley1990orlibrary}. JSPLib \citep{jsplib} and FJSPLib \citep{fjsplib} record the best-known lower and upper bounds of each instance, together with the solver and the hardware that produced them. JSPLib also publishes a standard 10-minute comparison of CP Optimizer, CP-SAT, OptalCP, and CPLEX on every instance. \citet{naderi2023} and \citet{yuraszeck2026} have compared constraint programming and mixed-integer programming solvers on thousands of shop instances, and \citet{dacol2022} have generated very large Taillard-type job shops to test how CP Optimizer and CP-SAT scale. Good practice for computational experiments requires public data, independent verification, and the reporting of all runs \citep{barr1995designing,hooker1995testing,gleixner2021miplib}. We follow these rules in this paper.

\subsection{Flexible job shops with transport robots}

In many plants, automated guided vehicles (AGVs) move jobs between machines, so machines and vehicles must be scheduled together \citep{bilge1995,hosseini2024}. The same integration of production and transport scheduling also arises along supply chains \citep{frazzon2018hybrid}. \citet{homayouni2021} have proposed a standard benchmark for this problem, the flexible job shop with transport (FJSPT). It adds a load/unload station, travel times, and two vehicles to the instances of \citet{brandimarte1993}, \citet{chambers1996}, and others. The small instances can now be solved exactly \citep{ham2020,homayouni2021,yao2025milp,yao2026cp,fjspagvmilp}. For the large instances, several heuristic methods have published results: late-acceptance hill climbing (LAHC) \citep{homayouni2021,burke2017}, a biased random-key genetic algorithm \citep{homayouni2023}, an improved genetic algorithm (IGA) \citep{meng2023}, a dual population genetic algorithm (DCGA) \citep{han2024dcga}, a tabu search \citep{berterottiere2024}, a variable neighborhood descent genetic algorithm (IVNDGA) \citep{jia2025ivndga}, a deep reinforcement learning method (HA-DQN) \citep{dong2025hadrl}, and an ant colony hybrid \citep{he2025aco}. \citet{meng2025cp} have proposed constraint programming models and a metaheuristic assisted by these models for this benchmark, and \citet{liu2026hgnn} have compared a learned policy with the best values of several of these methods. In this literature, the problem is also called the FJSP-AGV.

\subsection{LLMs for optimization and scheduling}

LLMs support optimization in two main ways. In the first, the LLM writes a model that a standard solver then solves, as in NL4Opt \citep{ramamonjison2023nl4opt}, OptiMUS \citep{ahmaditeshnizi2024optimus}, Chain-of-Experts \citep{xiao2024chainofexperts}, and ORLM \citep{huang2025orlm}. Alternatively, the LLM writes or improves heuristics inside an evolutionary loop, as in FunSearch \citep{romera2024funsearch}, OPRO \citep{yang2024opro}, EoH \citep{liu2024eoh}, ReEvo \citep{ye2024reevo}, AlphaEvolve \citep{novikov2025alphaevolve}, and the LLM4AD platform \citep{liu2024llm4ad}. Recent literature reviews can be found in \citet{liu2026algdesignsurvey} and \citet{daros2026llmcosurvey}.

In scheduling, LLMs have been used to generate schedules directly \citep{abgaryan2024llmsschedule,abgaryan2025starjob}, to design operators and dispatching rules for the flexible job shop \citep{liao2025onlineoperator,cao2025reflecsched}, and to coordinate multi-agent scheduling pipelines \citep{wang2025masc}. \citet{sanchezfernandez2026agentic} have tested whether LLM agents can write correct constraint programming models for six industrial scheduling problems. However, these studies mostly measure the correctness of models, or compare against simple baselines on small instances.

Our study differs from these papers in three ways. First, the agent writes a complete solution method, not only a model or a single operator. Second, the targets are large open instances with strong published bounds, most of them from commercial solvers. Third, the agent works in a general coding environment, without a special framework, an evolutionary loop, or task-specific training.

\section{Problem description}\label{sec:problem}

\subsection{Problem definitions}

We now formally introduce the three problems studied in this paper. In the job shop scheduling problem, a set $J$ of $n$ jobs must be processed on a set $M$ of $m$ machines. Each job $j\in J$ is a sequence of operations, and operation $o$ of job $j$ requires machine $\mu_{j,o}$ for $p_{j,o}$ time units without preemption. The operations of a job must follow this sequence, and a machine processes at most one operation at a time. The aim is to determine the start time $s_{j,o}$ of each operation, so that all constraints are satisfied and the makespan $C_{\max}=\max_{j,o}\{s_{j,o}+p_{j,o}\}$ is minimized. Let $L_j=\sum_{o}p_{j,o}$ be the length of job $j$. The value $\max_{j\in J}L_j$ is a lower bound on $C_{\max}$.

The flexible job shop scheduling problem extends this problem with a choice of machine. Operation $o$ of job $j$ can be processed on any machine $k$ of a set $M_{j,o}\subseteq M$, with a processing time $p_{j,o,k}$. The aim is to determine the machine and the start time of each operation, so that the makespan is minimized.

The FJSPT adds a load/unload station and a set $V$ of identical vehicles to the flexible job shop scheduling problem. Location 0 represents the station, and $t_{a,b}$ is the travel time from location $a$ to location $b$, with $a,b\in\{0\}\cup M$. All jobs and all vehicles are at location 0 at time zero. Each operation requires a loaded trip from the location of its predecessor to its machine, where location 0 is the predecessor of the first operation of a job. No trip is required if two consecutive operations of a job use the same machine. A trip departs after the completion of the preceding operation, and the operation starts after the arrival of the trip. A vehicle carries one job at a time, travels empty from each drop-off to its next pickup, and may wait. The buffers at the machines are unlimited, and the jobs do not return to location 0. The aim is to determine the machine of each operation, the vehicle of each trip, and all start times, so that the makespan is minimized. All instances of this paper have $|V|=2$. The three problems are NP-hard since they contain the job shop scheduling problem as a special case \citep{garey1976}.

\subsection{Lower bounds for the FJSPT}\label{sec:fjsptlb}

The best published lower bounds of the FJSPT \citep[Eqs.~26--27]{homayouni2021} are the longest path of a job through processing and loaded travel, and the total shortest loaded travel shared among the vehicles. They do not include empty travel. However, jobs never return to location 0, so the vehicles must make at least $n-|V|$ empty trips to this location. We now derive a lower bound based on this observation. We charge each empty trip at the length of a shortest path between its end points. This is conservative, because on the layouts that violate the triangle inequality, it is shorter than the direct travel time of the benchmark. Let $W_v$ be the total loaded and empty travel time of vehicle $v$ in a schedule, and let $q_v$ be the processing time of the operation that $v$ delivers last plus the shortest remaining processing and loaded travel time of its job. The last drop-off of $v$ occurs no earlier than $W_v$, so $C_{\max}\ge W_v+q_v$, and thus $|V'|\,C_{\max}\ge\sum_{v\in V'}(W_v+q_v)$ for the set $V'$ of vehicles used. In addition, let $\ell_k$ be the total processing time assigned to machine $k$, and let $h_k$ and $t_k$ be the smallest earliest start and the smallest remaining time after completion among the operations that are eligible on $k$, computed with the operation on $k$. If machine $k$ is used, then $C_{\max}\ge h_k+\ell_k+t_k$. Finally, $C_{\max}$ is at least the processing and loaded travel time along the path of each job. We minimize $C_{\max}$ subject to these inequalities over all machine assignments and all vehicle routes, without any other timing constraint. This yields a small integer program, in which each job is a path through its eligible machines, the empty trips balance the loaded trips at each location, and at most $|V|$ routes start at location 0. Its optimal value, rounded up to the time granularity of the instance, is a lower bound on $C_{\max}$. We refer to it as bound~(1).

The first of these inequalities yields a closed-form bound, bound~(2), that requires no solver. Let $\pi_a$ be the length of a shortest path from location $a$ to location 0. An empty trip from $a$ to $b$ takes at least $\pi_a-\pi_b$ time units. If we sum these differences over all empty trips, the terms of the loaded trips telescope along each job. Thus, the total empty time is at least $\sum_{j\in J}\pi_{f_j}-\sum_{v\in V}\pi_{e_v}$, where $f_j$ is the last machine of job $j$ and $e_v$ is the location of the last drop-off of vehicle $v$. Let $A_j$ be the minimum, over all machine sequences of job $j$, of its loaded travel time plus $\pi_{f_j}$. Let $\tau$ be the minimum, over all machines $k$, of the smallest processing time on $k$ minus $\pi_k$. Then $C_{\max}\ge\sum_{j\in J}A_j/|V|+\tau$.

When processing times dominate travel times, a second relaxation is more useful. We remove the capacity of the vehicles and the empty trips, and we keep each loaded trip as a delay between two consecutive operations of a job. The result is a flexible job shop with travel delays, which we solve by means of CP-SAT. Its optimal value is also a lower bound on $C_{\max}$, because every schedule of the FJSPT is feasible for this relaxation. We refer to it as relaxation~(3).

\section{Problem selection}\label{sec:selection}

We selected problems that are established, have tracked bounds, are open, are large and complex, and are checkable. They have been in the literature for several years with public instance files, a public table or a recent paper gives their best-known bounds and the solver that found them, commercial solvers are often far from the lower bound, and the feasibility of a schedule can be checked easily with independent code.

Table~\ref{tab:instances} lists the selected instances. The job shop \texttt{tai\_100\_1000\_1} has 100 jobs, 1{,}000 machines, and 100{,}000 operations \citep{dacol2022}. JSPLib lists it as open, with recent bounds from OptalCP, CP Optimizer, and Hexaly \citep{jsplib,yuraszeck2026}. From the Behnke--Geiger flexible job shops \citep{behnke2012}, we have selected the open instances \texttt{behnke16} and \texttt{behnke14}, and the closed instance \texttt{behnke31} as a control, since a good method should find its known optimum. For the flexible job shop with transport robots, we have used the benchmark of \citet{homayouni2021}, the standard benchmark for this problem, on which several later methods report results \citep{meng2023,berterottiere2024,dong2025hadrl}. The problem adds a second, interacting resource (the vehicles) to the flexible job shop, which makes it a good representative of integrated production and logistics scheduling \citep{frazzon2018hybrid}. From this benchmark, we have selected the ten MKT instances and the 21 Chambers--Barnes instances. They are the largest sets, and published heuristics were far from the published lower bounds. However, \citet{meng2025cp} report 14 of the 21 Chambers--Barnes instances as solved to optimality with CP Optimizer (Section~\ref{sec:res-fjspt}). The smaller sets are mostly solved to optimality, so we have used them only to test generalization. In Table~\ref{tab:instances}, LB and UB are the best-known lower and upper bounds before this work. For the transport instances, LB is the best published lower bound \citep{homayouni2021}. The asterisk marks the control instance.

\begin{table}[!htbp]
\caption{Selected instances}
\label{tab:instances}
\footnotesize
\begin{tabular*}{\textwidth}{@{\extracolsep{\fill}}l l r r r r l@{}}
\toprule
Problem & Instance(s) & Jobs & Machines & Operations & LB / UB & Status \\
\midrule
Job shop & \texttt{tai\_100\_1000\_1} & 100 & 1{,}000 & 100{,}000 & 522{,}298 / 533{,}080 & open \\
Flexible job shop & \texttt{behnke16} & 100 & 20 & 500 & 391 / 412 & open \\
 & \texttt{behnke14} & 50 & 20 & 250 & 221 / 230 & open \\
 & \texttt{behnke31} & 50 & 40 & 250 & 226 / 226 & closed$^{\ast}$ \\
Flexible job shop, 2 AGVs & MKT01--MKT10 & 10--20 & 4--15 & 55--240 & e.g.\ MKT10: 386.8 / 525 & open \\
 & 21 Chambers--Barnes & 10--15 & 11--18 & 100--225 & see Table~\ref{tab:cb} & 7 open \\
\bottomrule
\end{tabular*}
\end{table}

\section{Method}\label{sec:method}

\subsection{Framework and technology stack}

\begin{figure}[!htbp]
\centering
\begin{tikzpicture}[scale=0.84, transform shape,
  box/.style={draw, rounded corners=2pt, align=center, font=\small, minimum height=1.0cm, text width=3.0cm, fill=#1},
  box/.default=white,
  arr/.style={-{Stealth[length=2mm]}, thick}]
\node[box=blue!6] (task) {\textbf{Task}\\instance file, rules,\\target, time budget};
\node[box=orange!10, right=1.0cm of task, text width=4.6cm, minimum height=2.2cm] (agent) {\textbf{LLM coding agent}\\(Claude Code)\\[2pt] analyze instance $\cdot$ write code\\compile and run $\cdot$ read output\\change the method};
\node[box=gray!8, above=0.5cm of agent, text width=4.6cm] (tools) {open-source tools: Python, C/C++,\\NumPy, OR-Tools CP-SAT};
\node[box=green!8, right=1.0cm of agent] (out) {\textbf{Output}\\schedule (mode A)\\or solver (mode B)};
\node[box=red!6, right=0.8cm of out] (check) {\textbf{Independent}\\\textbf{checker}\\(no shared code)};
\node[box=gray!8, below=0.55cm of agent, text width=4.6cm] (human) {human input (assisted runs only):\\goals, literature, compute};
\draw[arr] (task) -- (agent);
\draw[arr] (tools) -- (agent);
\draw[arr] (agent) -- (out);
\draw[arr] (out) -- (check);
\draw[arr, dashed] (human) -- (agent);
\draw[arr] (agent.north east) ++(0,-0.2) to[out=20,in=-20,looseness=4] node[right, font=\scriptsize]{loop} ([yshift=-0.6cm]agent.north east);
\end{tikzpicture}
\caption{Framework of the study}
\label{fig:framework}
\end{figure}
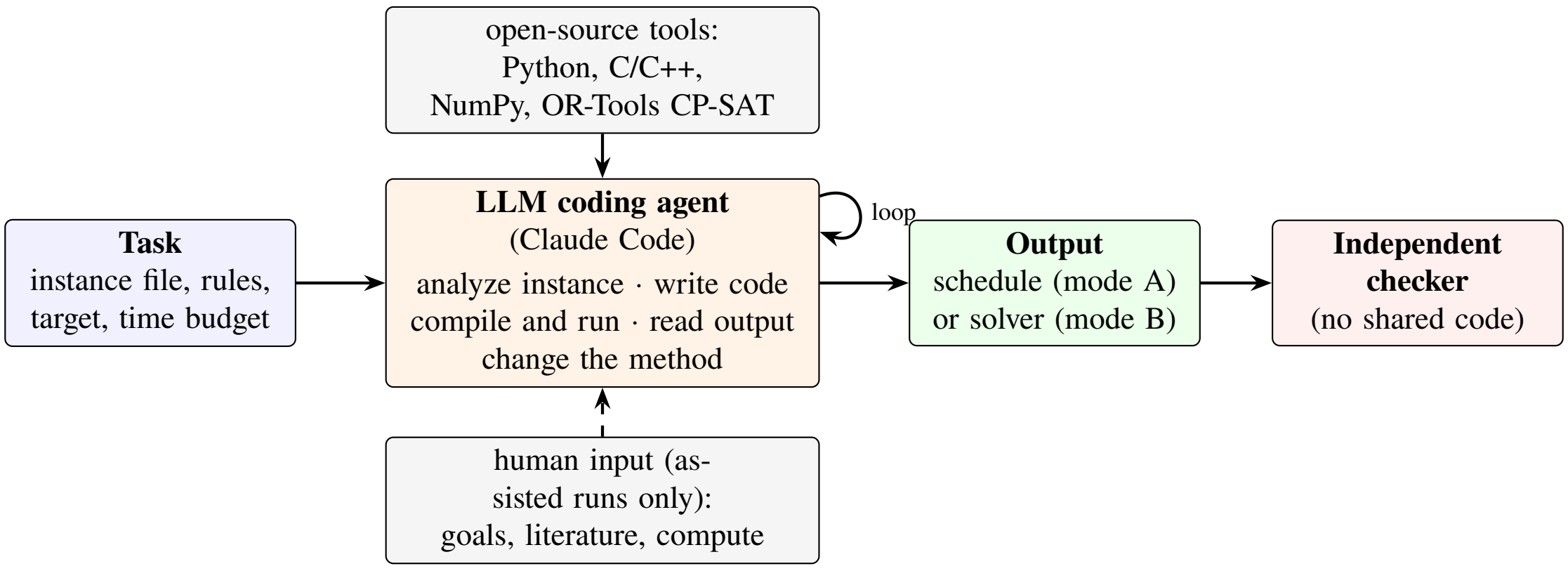

Figure~\ref{fig:framework} shows the framework of our study. In autonomous runs, the task is the only input. In assisted runs, we have also provided goals and papers, and we have managed the compute resources (dashed arrow). An independent checker verifies every output schedule. The agent is Claude Code \citep{claudecode}, a command-line coding agent. It reads files, writes code in any language, compiles and runs it, reads the output, and repeats this loop until the time budget ends. Table~\ref{tab:stack} lists the models, the hardware, and the software.

\begin{table}[!htbp]
\caption{Models, hardware, and software}
\label{tab:stack}
\footnotesize
\begin{tabular*}{\textwidth}{@{\extracolsep{\fill}}p{0.22\textwidth}p{0.74\textwidth}@{}}
\toprule
Item & Value \\
\midrule
Agent & Claude Code \citep{claudecode}, version 2.1, in all sessions \\
Models of record & \texttt{claude-opus-5} (Claude Opus~5) and \texttt{claude-sonnet-5} (Claude Sonnet~5): job shop and flexible job shop, accessed 5--10 August 2026; \texttt{claude-opus-5-5} (Claude Opus~5.5): transport robots, accessed 25--26 September 2026, and job shop in mode~C, 19--20 October 2026; \texttt{claude-fable-5-1} (Claude Fable~5.1): flexible job shop in mode~C, 19--20 October 2026. Public release \citep{anthropicreleasenotes}: Claude Sonnet~5, 30 June 2026; Claude Opus~5, 24 July 2026; Claude Fable~5.1, 1 September 2026; Claude Opus~5.5, 22 September 2026. Default decoding settings of Claude Code; the provider does not expose sampling seeds. \\
Autonomous runs & Harbor task framework \citep{harbor}; one Docker container per run; quota of 4 CPUs (three runs: 8 CPUs) and 4\,GB of memory, although the container reported the 32 hardware threads of the host; budget 18{,}000\,s (5\,h) \\
Container software & Python~3.12.3, GCC~13.3.0, NumPy~2.4.6; OR-Tools only if the agent installed it \\
Workstation & AMD Ryzen~9 5950X (3.4\,GHz), 16 cores / 32 threads, 125\,GB RAM; job shop in mode~C: 4-core ARM cloud machine (Neoverse-N1) \\
Workstation software & Python~3.11, NumPy, GCC (C/C++), OR-Tools CP-SAT~9.15 \\
Commercial solvers & none used by the agent; CP Optimizer~22.2 and Hexaly~15.0 (academic licenses) for our baseline runs \\
\bottomrule
\end{tabular*}
\end{table}

\subsection{Autonomous runs and assisted runs}

Our main goal is to verify whether the agent is able to obtain good solutions without expert knowledge. We therefore distinguish two types of runs.

In an \emph{autonomous run}, the agent was given one task: improve the makespan of one instance. The task text described the instance format, the rules, the output format, and a target. For \texttt{tai\_100\_1000\_1}, the target was a makespan strictly below 540{,}000. For \texttt{behnke16}, the text asked for a heuristic, stated the best-known upper bound of 412, and set the target 411. The agent then worked without human input in a new container, with a fixed amount of time and of CPU cores, and it selected its own method. An independent test checked the schedule at the end. One job-shop session also received code written by an earlier session and the best value then known.

In an \emph{assisted run}, we worked with the agent in an interactive session on the workstation. We gave goals and feedback, supplied papers, and managed the compute resources (for example, how many cores a solver could use), but we did not write any solver code. Table~\ref{tab:runs} shows the type of each run.

\subsection{Modes of use of the agent}\label{sec:modes}

We have tested three modes of use. In mode~A, the agent itself is the algorithm, and it must deliver a good schedule for one instance. It may write code that is used only once, and it does not have to leave a reusable program. If a new schedule is needed, the agent is run again. In mode~B, the agent must deliver a solver that the user can compile and run again, on the same instance or on similar ones. The cost then shifts from LLM tokens to CPU time. In mode~C, the agent starts from the outcome of modes A and B, namely the solvers and best schedules of this paper, and runs its own improvement campaign with the goal of closing as many instances as possible. It may change the solvers, write new code, and select the experiments to run. We have run mode~C once on the job shop and once on the flexible job shop, after the other experiments, each as a single, mostly autonomous session that used open-source tools only. We set the goal and stopped the session, but we suggested no method and wrote no code. Section~\ref{sec:disc} compares the modes.

\subsection{Verification}\label{sec:verif}

We have not relied on the checkers written by the agents. For each problem family, a separate checker, which shares no code with the solvers, verifies every schedule. The FJSPT checker uses exact fractions, because some layouts use half time units. We have also tested the checkers on deliberately corrupted schedules, and they rejected all of them. Each instance file was compared with the original distribution and with independent copies. A single script in the public repository re-checks every result of this paper. We take the lower bounds of the libraries and the literature as published, but not the single-solver verdicts of our baseline runs. An instance is closed only when a certified schedule reaches such a bound or a bound of Section~\ref{sec:fjsptlb} proved by two solvers.

\subsection{Contamination controls}\label{sec:contamination}

The instances are public, and they may be part of the training data of the models. However, no schedule with the makespans reported here existed before this work, so a new best-known value cannot be reproduced from memory, and every value is certified by the independent checker. The containers of the autonomous runs had public network access, but no session used a web search or fetch tool, and no session downloaded benchmark data, bounds, or published solutions.  The autonomous runs received only the instance and the task text, and for \texttt{behnke16} a small example instance and schedule. One job-shop session also received earlier code (Section~\ref{sec:method}). The target values of the job-shop runs came from the first interactive session (Section~\ref{sec:reporting}). The assisted transport-robot session started from one message that pasted the paper of \citet{meng2023} and asked the agent to choose the instance with the best chance of improvement and to find tailored solutions rather than to rely on computation. During the session, we supplied the further papers that it requested, including published schedules. It used these schedules only to validate its model and data (Section~\ref{sec:how}). No commercial solver was available to any agent.

\subsection{Reporting and repetitions}\label{sec:reporting}

Table~\ref{tab:runs} lists every autonomous run and every interactive agent session on the instances of this paper, as they were recorded, including the runs that missed the target, failed, or stopped early. No run has been removed on the basis of its result. Runs on other problems, and on instances that we generated to calibrate the tasks, are not included. The instance \texttt{tai\_100\_1000\_1} was chosen in the first interactive session, after a greedy pilot run, and the target values of the autonomous job-shop runs came from that session. For the assisted Behnke session, no transcript was kept. The transport-robot study was one long assisted session in which the agent made 380 solver runs. For it, we report the best schedules, a fixed protocol of 100 runs (Table~\ref{tab:mkt}), and the controlled comparisons of moves (Sections~\ref{sec:how} and~\ref{sec:res-fjspt}).

Because the provider does not expose sampling seeds, we measure the run-to-run variation of the agents by repeated sessions with the same task. On \texttt{tai\_100\_1000\_1}, the six Claude Opus~5 sessions, one of which received earlier code (Section~\ref{sec:method}), ended between 527{,}000 and 533{,}386 (median 531{,}403), and the nine Claude Sonnet~5 sessions that delivered a schedule between 535{,}930 and 554{,}332 (median 541{,}641). On \texttt{behnke16}, the twelve autonomous sessions ended between 404 and 479 (median 418). The solvers written by the agents are re-run with explicit seeds, and we report the best value together with the range or the median of the runs (Tables~\ref{tab:jsspgen} and~\ref{tab:mkt}).

\section{Computational experiments}\label{sec:results}

All computations on our workstation were executed on an AMD Ryzen~9 5950X (Table~\ref{tab:stack}).

Table~\ref{tab:runs} lists all agent runs on the instances of this paper and all re-runs of the solvers written by the agents. In this table, ``Auto'' denotes an autonomous run and ``Assist'' an assisted run. The column ``Mode'' shows whether the agent is the solver and the output is a schedule (A), whether the output is a reusable solver (B), or whether the agent improves the solvers of this paper in its own campaign (C). The number of independent runs is given in the column ``Runs''. The 12 runs of the 527{,}000 solver ran in parallel with one core each, and the 480 runs of the corridor solver are the workers of its 30 runs with 16 cores. Best values are certified by the independent checker.

\begin{table}[!htbp]
\caption{All agent runs and re-runs of the agent solvers}
\label{tab:runs}
\scriptsize
\setlength{\tabcolsep}{3pt}
\begin{tabular*}{\textwidth}{@{\extracolsep{\fill}}l l l l r l l@{}}
\toprule
Instance(s) & Type & Mode & Model & Runs & Result (makespan) & Human input \\
\midrule
\texttt{tai\_100\_1000\_1} & Auto & A & Opus 5 / Sonnet 5 & 16 & best 527{,}000; 15 feasible & none\tabnotemark{a} \\
\texttt{tai\_100\_1000\_1} & Assist & A & Opus 5 & 1 & 530{,}404 & goals, review \\
\texttt{tai\_100\_1000\_1} & Assist & B & --\tabnotemark{b} & 12 & 525{,}764 & compute \\
\texttt{behnke16} & Auto & A & Opus 5 & 5 & 404, 411, 415, 417, 419 & none \\
\texttt{behnke16} & Auto & A & Opus 5 & 5 & 423, 430, 435, 450, 479\tabnotemark{c} & none \\
\texttt{behnke14/16/31} & Auto & A & Opus 5 & 2 & 234--235 / 415--416 / 226--229 & none \\
\texttt{behnke14/16/31} & Assist & B & Opus 5 & 1 & 232 / 406 / 226 & goals, review \\
MKT01--10, 21 Chambers--Barnes & Assist & B & Opus 5.5 & 1 & Tables~\ref{tab:mkt} and~\ref{tab:cb} & goal, papers, review \\
\texttt{tai\_100\_1000\_1..10} & Assist & C & Opus 5.5 & 1 & 4 of 10 optimal; Table~\ref{tab:modec} & goal, choice of instances \\
Behnke 1--60 & Assist & C & Fable 5.1 & 1 & 6 newly closed; Table~\ref{tab:modec} & goal, tools \\
\midrule
\multicolumn{7}{@{}l}{\emph{Re-runs of the agent solvers, without change (no agent session)}} \\
\texttt{tai\_100\_1000\_1..10} & -- & B & corridor solver & 30 + 128 + 480 & 10/10 below best known; 2 optimal & none \\
Behnke 1--60 & -- & B & reusable tabu & 60 + 32 & 11 below best known & none \\
MKT01--10 & -- & B & transport solver & 100 & 61/100 below reference value & none \\
SFJST/MFJST 01--10 & -- & B & transport solver & 100 & all small optima & none \\
\bottomrule
\end{tabular*}
\tablenote{\tabnotemark{a}One session (530{,}881) also received code written by an earlier session and the best value then known (Section~\ref{sec:method}). \tabnotemark{b}The solver of the 527{,}000 session, recovered from its log and run on the workstation. \tabnotemark{c}Three of these runs stopped after 4--7\% of the budget, and one stopped after 462\,s because of an API error. Only their reports and final values were kept.}
\end{table}

\subsection{How the agents found their solutions}\label{sec:how}

The agents have not applied one fixed algorithm. In the successful runs, the agent first analyzed the instance and then adapted the formulation or the method to its structure. This subsection describes the main ideas, and Sections~\ref{sec:res-jssp}--\ref{sec:res-newbest} report the results.

\subsubsection{Change of formulation for the job shop}

\begin{figure}[!htbp]
  \centering
  \includegraphics[width=0.75\textwidth]{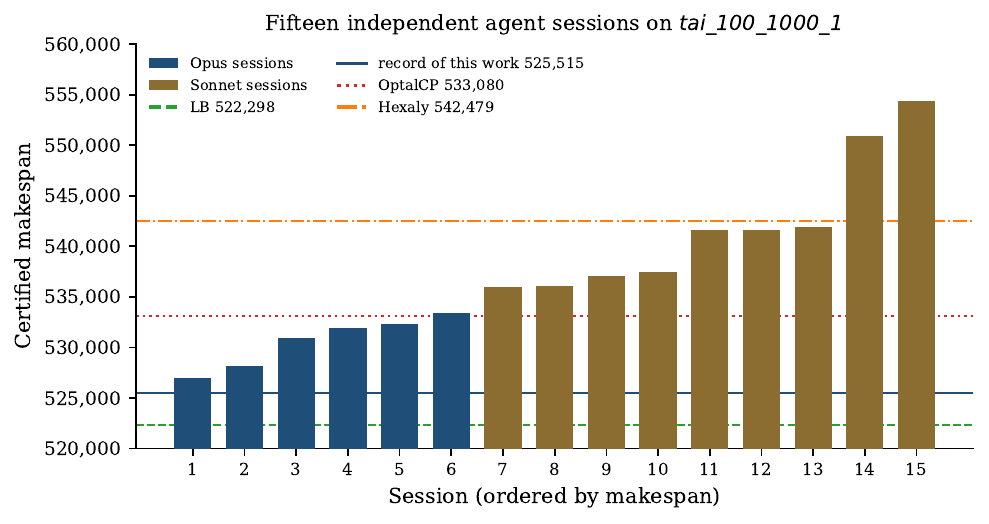}
  \caption{Final makespans of the autonomous sessions on \texttt{tai\_100\_1000\_1}}
  \label{fig:sessions}
\end{figure}

\begin{figure}[!htbp]
  \centering
  \includegraphics[width=0.495\textwidth]{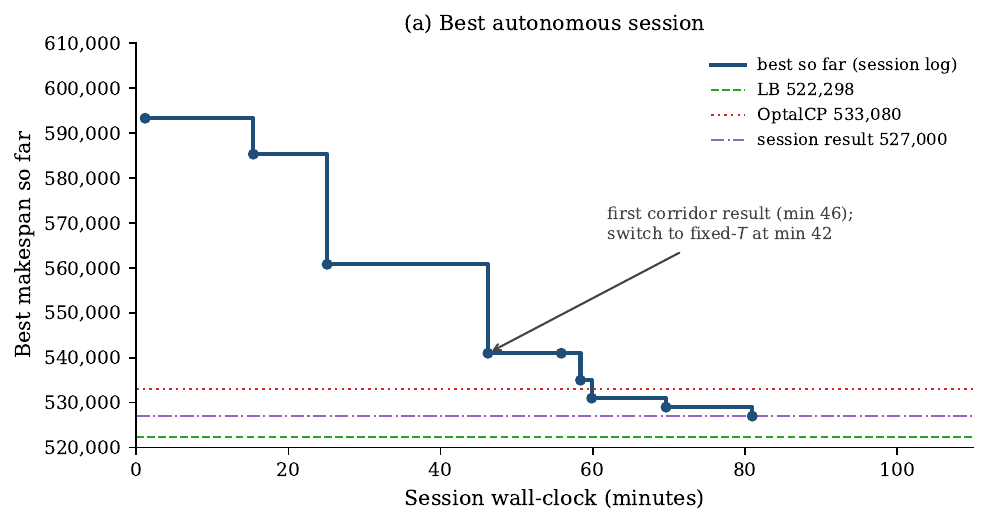}\hfill
  \includegraphics[width=0.495\textwidth]{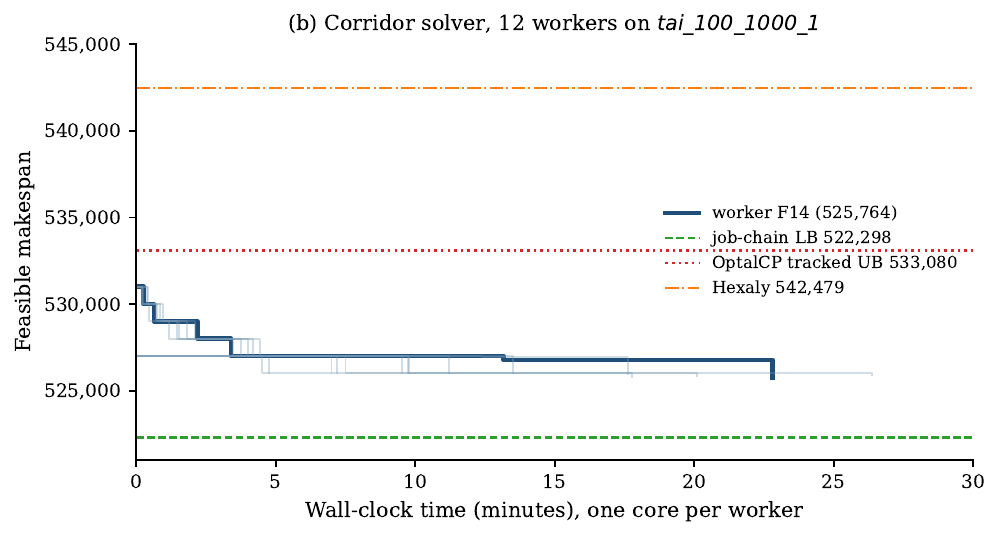}
  \caption{Best makespan on \texttt{tai\_100\_1000\_1} over time}
  \label{fig:timeline}
\end{figure}

In \texttt{tai\_100\_1000\_1}, each job visits 1{,}000 machines. The machines are idle about 90\% of the time, while each job is busy for about 91--99\% of the time. The longest job alone needs 522{,}298 time units, which is the lower bound. The best of the autonomous sessions (Figure~\ref{fig:sessions}) started from a greedy schedule of 593{,}321 and wrote four solvers that minimized the makespan, and the best of them reached 592{,}865. At minute~42, the agent changed the formulation. It fixed a target makespan $T$ and searched for a schedule with no machine conflict and $C_{\max}\le T$. To this end, the agent restricted each operation of job~$j$ to a corridor
\[
\mathrm{pref}_j(o) \;\le\; s_{j,o} \;\le\; \mathrm{pref}_j(o) + (T-L_j),
\]
where $L_j$ is the length of job~$j$ and $\mathrm{pref}_j(o)$ is the start of operation~$o$ when the job runs without idle time. Inside these corridors, $C_{\max}\le T$ always holds. The corridors alone do not ensure the job order, because two consecutive operations of a job may overlap while each of them lies inside its own corridor. The search therefore starts from the schedule $s_{j,o}=\mathrm{pref}_j(o)+(T-L_j)\,o/m$, in which the delay of an operation does not decrease along its job, and it maintains the job order at every iteration. When it shifts an operation, it also shifts the successors or the predecessors of this operation on the job as far as required, and it rejects the shift if an operation would leave its corridor or if more than 60 operations would move. Thus, only machine conflicts remain. A min-conflicts search with random perturbations removes them, and each time it finds a feasible schedule, the agent lowers $T$. At each iteration, the search selects two overlapping operations on a machine, evaluates the four shifts that separate them, and applies the shift that yields the smallest total overlap, or a random feasible shift with a probability of 8\%. After 200{,}000 iterations without improvement, it shifts random operations of machines with overlaps. After a feasible schedule, $T$ decreases by 1{,}000, and after 40 million iterations without improvement, the search restarts from the best schedule with a random target below its makespan. The first corridor prototype, started at $T=545{,}000$, found a schedule of 540{,}999 after about 60\,s, at minute~46 of the session. The agent then wrote a faster C version (\texttt{csp3.c}) and ran 28 worker processes (16 of \texttt{csp3} and 12 of a variant) on the 4 CPUs of the container. A \texttt{csp3} worker found 527{,}000 at about minute~77, and the agent saved it at minute~81 of a session that lasted 101.6 minutes. Figure~\ref{fig:timeline}a shows this session.

A second session searched over non-decreasing delay profiles of the jobs and found 528{,}142. An assisted campaign combined tabu search with sequence fixing by CP-SAT and found 530{,}404. The corridor method yields better values than both.

\subsubsection{Stage structure of the flexible job shop}

The best session on \texttt{behnke16} analyzed the eligibility sets and found that the instance is a three-stage flexible flow shop: operation~0 of each job runs on machines 1--4, operations 1--3 on machines 5--16, and operation~4 on machines 17--20. The last stage is the bottleneck. The agent then wrote a parallel tabu search with ruin-and-recreate moves. After each move, a fast procedure re-optimizes the machine assignment of the first and the last stage. With the rest of the schedule fixed, each of these stages reduces to a partition problem. A secondary criterion counts the operations near the critical path and breaks ties among schedules of equal makespan. The session improved a constructive solution of 448 down to 404. It ran about 31 worker threads in a container with a quota of 4 CPUs. The other four sessions of this batch used tabu search, simulated annealing, CP-SAT, and large neighborhood search, and ended between 411 and 419. One session used CP-SAT with a warm start and proved that 226 is optimal for \texttt{behnke31}.

\subsubsection{Routing structure of the problem with transport robots}

In the assisted FJSPT session, we gave the agent the benchmark, one recent paper with published results \citep{meng2023}, and the goal of improving the best-known values. We also asked it to exploit the structure of the problem rather than rely on more computation, and during the session we supplied the recent papers that it requested. The agent then executed the following steps:
\begin{enumerate}
\item It used MKT10, the largest MKT instance, as the main development instance.
\item It replayed the published schedules of two competing papers \citep{meng2023,berterottiere2024} with an independent checker. All of them reproduced the published makespans exactly, which confirmed that the model and the data match the literature.
\item It analyzed good schedules and found that both vehicles are almost never idle and that most arcs of the critical path connect consecutive tasks of the same vehicle. In other words, the problem is mainly a two-vehicle pickup-and-delivery routing problem.
\item It wrote a simulated annealing solver with moves borrowed from vehicle routing: relocation, exchange, 2-opt$^\ast$ \citep{potvin1995}, and ruin-and-recreate \citep{schrimpf2000}.
\item It ran controlled experiments, kept the moves that improved the results, and removed the others.
\end{enumerate}

\subsection{Job shop}\label{sec:res-jssp}

\paragraph{Autonomous sessions (mode A)} Sixteen autonomous sessions worked on \texttt{tai\_100\_1000\_1}: six with Claude Opus~5 and ten with Claude Sonnet~5. Fifteen of them delivered a feasible schedule, and one stopped after 18 minutes without a schedule. Figure~\ref{fig:sessions} shows the final makespan of each session, numbered in order of makespan. The horizontal lines show the job-chain lower bound \citep{jsplib}, the best value of this work, the tracked OptalCP bound \citep{jsplib,optalcp}, and the Hexaly bound \citep{yuraszeck2026}. Five of the six Opus sessions finished below the tracked OptalCP bound of 533{,}080. All six Opus sessions ended below the nine Sonnet sessions that delivered a schedule (at most 533{,}386 against at least 535{,}930). All 15 schedules are below the 10-minute CP-SAT value (554{,}381), and 13 are below the 10-minute CP Optimizer value (545{,}805). The best session reached 527{,}000 with 4 CPUs, in 101.6 minutes and with an API cost of US\$10.07.

\paragraph{Agent solver (mode B)} The best session left its solver, \texttt{csp3.c}, in the container. We have recovered the source from the session log and have run it with one core per worker. Twelve workers ran for 1{,}800\,s each. All twelve finished below the 527{,}000 of the session, and five below 526{,}000. The best worker (F14 in Figure~\ref{fig:timeline}b) reached 525{,}764 after 23 minutes. Six workers started with the target 531{,}000 and six with 527{,}000, the result of the session. The best worker was one of the first six. The logs do not record the host, but the iteration rates match one core per worker. Within the session, the solver did not reach 526{,}000. This is due to the fact that its 28 processes shared 4 CPUs, which reduced the iteration rate of each process to about one tenth of that of a process with its own core. On a single core, the same solver reaches 533{,}000, the tracked OptalCP value, in 2.7--5.1\,s when it starts with this target, and in 7.8--16.2\,s when it starts at 535{,}000. With the 16-core budget of our CP Optimizer runs (16 workers with their own cores and different seeds for 1{,}800\,s, three such runs), the same solver reached 525{,}515, the best value of this work.

\begin{table}[!htbp]
\caption{Solver and agent results on \texttt{tai\_100\_1000\_1}}
\label{tab:jssp}
\footnotesize
\begin{tabular*}{\textwidth}{@{\extracolsep{\fill}}l l r r@{}}
\toprule
Method & Budget & Makespan & Gap (\%) \\
\midrule
CP-SAT 9.15 \citep{jsplib} & 10 min, laptop & 554{,}381 & 6.14 \\
CP Optimizer 22.1 \citep{jsplib} & 10 min, laptop & 545{,}805 & 4.50 \\
OptalCP 2026.4 \citep{jsplib} & 10 min, laptop & 540{,}099 & 3.41 \\
CPLEX 22.1 \citep{jsplib} & 10 min, laptop & none & -- \\
CP Optimizer 12.10 \citep{dacol2022} & 6 h, 1 core & 542{,}768 & 3.92 \\
Hexaly 14.0 \citep{yuraszeck2026} & 600\,s, 12 threads & 542{,}479 & 3.86 \\
OptalCP 2026.2 \citep{jsplib,optalcpbench} & 180\,min, 16 workers & 533{,}080 & 2.06 \\
\midrule
CP Optimizer 22.2 (ours) & 10 min, 1 core & 551{,}319 & 5.56 \\
CP Optimizer 22.2 (ours) & 30 min, 1 core & 546{,}035 & 4.54 \\
CP Optimizer 22.2 (ours) & 30 min, 16 cores & 543{,}034 & 3.97 \\
Hexaly 15.0 (ours) & 10 min, 1 core & 541{,}994 & 3.77 \\
Hexaly 15.0 (ours) & 30 min, 1 core & 541{,}764 & 3.73 \\
Hexaly 15.0 (ours) & 30 min, 16 cores & 539{,}566 & 3.31 \\
\midrule
Agent session (A) & 101.6 min, 4 CPUs & 527{,}000 & 0.90 \\
Agent solver (B) & 10 min, 1 core & 526{,}999 & 0.90 \\
Agent solver (B) & 30 min, 1 core & 525{,}999 & 0.71 \\
Agent solver (B) & 30 min, 12 $\times$ 1 core & 525{,}764 & 0.66 \\
Agent solver (B) & 30 min, 16 cores & \textbf{525{,}515} & \textbf{0.62} \\
\bottomrule
\end{tabular*}
\end{table}

\paragraph{Comparison with solvers} Table~\ref{tab:jssp} compares these results with the published values of commercial and open-source solvers.  The gap is relative to the lower bound of 522{,}298. The JSPLib values come from a standard 10-minute run on a 4-core laptop processor (Intel Core i7-1185G7). JSPLib does not state the number of threads. The tracked OptalCP value comes from a run of 180 minutes. JSPLib lists the processor of our workstation (Ryzen~9 5950X) for this run, but the published run record of OptalCP gives version 2026.2, 16 workers, and an AMD EPYC 9454P \citep{optalcpbench}. The hardware, versions, and thread counts of the other published rows are those stated by the sources. \citet{dacol2022} used CP Optimizer~12.10 on an AMD EPYC 7551P, and \citet{yuraszeck2026} used Hexaly~14.0 on an Intel Core i7-10750H with all 12 threads. Because these published values come from different hardware and budgets, we have also run CP Optimizer~22.2 and Hexaly~15.0 ourselves on our workstation, with the same budgets as the agent solver (Tables~\ref{tab:jssp}, \ref{tab:jsspgen}, and~\ref{tab:compare}). These rows are labeled ``ours''. The 1-core and 16-core rows of our runs are the best of three runs, and a 16-core run of the agent solver consists of 16 single-core workers with different seeds. For both, we have used the standard job-shop model of the vendor with default parameters. Hexaly reached 541{,}764 with one core and 539{,}566 with 16 cores after 30 minutes, better than its published value of 542{,}479 and than CP Optimizer, but 2.7\% above the agent solver with 16 cores. Our Hexaly license allowed one process at a time, so each Hexaly run had the workstation to itself.

\begin{table}[!htbp]
\caption{Generalization of the corridor solver to the ten $100\times1000$ instances}
\label{tab:jsspgen}
\scriptsize
\setlength{\tabcolsep}{1.5pt}
\begin{tabular*}{\textwidth}{@{\extracolsep{\fill}}l r r r r r r r r r r r l@{}}
\toprule
 & & & & \multicolumn{2}{c}{CP Optimizer, 1 core} & \multicolumn{2}{c}{Agent, block 1} & Block 2 & Agent & & & \\
\cmidrule(lr){5-6}\cmidrule(lr){7-8}
Instance & LB & Best known & OptalCP 10 min & 10 min & 30 min & 10 min & 30 min & 16 seeds & 16 cores & Best & Gap (\%) & Status \\
\midrule
\texttt{\_1} & 522{,}298 & 533{,}080 & 540{,}099 & 551{,}319 & 546{,}035 & 526{,}999 & 525{,}999 & 526{,}000 & 525{,}515 & \textbf{525{,}515} & 0.62 &  \\
\texttt{\_2} & 530{,}375 & 538{,}067 & 544{,}180 & 551{,}476 & 550{,}373 & 532{,}000 & 530{,}375 & -- & 530{,}375 & \textbf{530{,}375} & 0.00 & optimal \\
\texttt{\_3} & 530{,}560 & 538{,}757 & 545{,}590 & 551{,}844 & 550{,}391 & 530{,}560 & 530{,}560 & -- & 530{,}560 & \textbf{530{,}560} & 0.00 & optimal \\
\texttt{\_4} & 527{,}101 & 534{,}746 & 541{,}326 & 551{,}836 & 546{,}753 & 529{,}999 & 529{,}075 & 529{,}000 & 528{,}156 & \textbf{528{,}156} & 0.20 &  \\
\texttt{\_5} & 517{,}728 & 529{,}580 & 537{,}153 & 552{,}148 & 547{,}705 & 527{,}000 & 527{,}000 & 526{,}333 & 525{,}865 & \textbf{525{,}865} & 1.57 &  \\
\texttt{\_6} & 522{,}907 & 534{,}969 & 542{,}004 & 549{,}319 & 547{,}290 & 526{,}000 & 525{,}995 & 525{,}000 & 524{,}812 & \textbf{524{,}812} & 0.36 &  \\
\texttt{\_7} & 522{,}537 & 534{,}974 & 541{,}161 & 550{,}879 & 547{,}556 & 527{,}000 & 527{,}000 & 527{,}805 & 527{,}000 & \textbf{527{,}000} & 0.85 &  \\
\texttt{\_8} & 526{,}428 & 535{,}757 & 543{,}074 & 551{,}819 & 549{,}836 & 530{,}000 & 529{,}282 & 528{,}817 & 528{,}473 & \textbf{528{,}473} & 0.39 &  \\
\texttt{\_9} & 528{,}097 & 536{,}993 & 544{,}065 & 552{,}682 & 552{,}462 & 532{,}000 & 531{,}207 & 531{,}000 & 529{,}470 & \textbf{529{,}470} & 0.26 &  \\
\texttt{\_10} & 521{,}766 & 529{,}918 & 538{,}445 & 551{,}704 & 545{,}227 & 527{,}000 & 526{,}855 & 527{,}000 & 526{,}000 & \textbf{526{,}000} & 0.81 &  \\
\bottomrule
\end{tabular*}
\end{table}

\paragraph{Generalization to similar instances} \texttt{tai\_100\_1000\_1} belongs to a family of ten instances of the same size and type, all of them open in JSPLib. We have run the corridor solver just described on all ten without any change, with three runs per instance with different seeds, one core, and 1{,}800\,s per run. Table~\ref{tab:jsspgen} shows the results. In Table~\ref{tab:jsspgen}, LB is the length of the longest job, ``Best known'' the JSPLib upper bound before this work, and ``OptalCP 10 min'' the standard JSPLib run \citep{jsplib}. CP Optimizer~22.2 and block~1 of the agent solver are the best of three seeds with one core after 10 and 30 minutes. The agent solver starts at $T = \lceil 1.02\,\mathrm{LB}/1000\rceil\cdot1000$. Block~2 and the 16-core runs are described below. In the column ``Best'', bold values are below the best-known value.

The solver found a schedule below the best-known upper bound on every instance, already after 10 minutes, far below the 10-minute values of OptalCP, CP Optimizer, and CP-SAT. With one core and 30 minutes, CP Optimizer~22.2 ended 3.7--5.8\% and Hexaly~15.0 2.7--4.4\% above the lower bound, and with 16 cores 2.6--5.3\% and 2.4--4.2\%. Both stayed above the old best-known value on every instance, whereas the agent solver ended 0.0--1.8\% and 0.0--1.6\% above the lower bound. Over the three runs and the ten instances, the mean gap with one core is 0.73\% for the agent solver, 4.49\% for CP Optimizer, and 3.60\% for Hexaly, and with 16 cores 0.53\%, 3.79\%, and 3.30\%. At both budgets, the worst run of the agent solver is below the best run of both commercial solvers on all ten instances. On \texttt{tai\_100\_1000\_2} and \texttt{\_3}, the solver found a makespan equal to the lower bound, so these two instances are now closed. Two of the three seeds found each optimum, that of \texttt{\_3} in less than 10 minutes. A second block of sixteen further runs of 1{,}800\,s with two processes per core improved five instances further but closed none, because each process had only half a core. Finally, with the 16-core budget of CP Optimizer (three runs of 16 single-core workers with different seeds), the solver found the optima of \texttt{\_2} and \texttt{\_3} in all runs and improved seven instances further. On \texttt{\_7}, all three runs ended at 527{,}000.

\paragraph{Improvement campaign (mode~C)} Starting from the corridor solver and the schedules above, one Claude Opus~5.5 session worked for 24\,h on a 4-core ARM cloud machine (about 84 core-hours of solver runs), with the goal of closing the ten instances. It made three changes to the solver. First, it started the search directly at the lower bound instead of descending in steps of 1{,}000. Second, it narrowed the time window of every operation by constraint propagation at that target. Third, it reduced the size of the random perturbations of the local search. With these changes, it found schedules whose makespan equals the lower bound on \texttt{\_4} (527{,}101) and \texttt{\_9} (528{,}097), so four of the ten instances are now closed, and it improved \texttt{\_8} to 527{,}280. Among the 191 logged runs of this campaign, one of the 34 runs with the target at the lower bound succeeded on \texttt{\_4}, one of 34 on \texttt{\_9}, none of 58 on \texttt{\_8}, and none of 23 on \texttt{\_6}. It also found that our schedules contain small idle gaps, and starting every operation as early as possible while keeping the order on each machine has lowered the other five values by 126--229 units. The session re-found the optimum of \texttt{\_2} in 47\,s on one core. Table~\ref{tab:modec} in Appendix~\ref{app:newbest} lists all values.

\subsection{Flexible job shop}\label{sec:res-fjsp}

Table~\ref{tab:behnke} lists all runs on the three Behnke instances (bold: below the best-known upper bound of FJSPLib). Of ten autonomous sessions on \texttt{behnke16} with a target of 411, the best reached 404 and a second one 411, both below the best-known 412. Three sessions stopped by themselves after 4--7\% of the budget, and one after 462\,s because of an API error. Of two further sessions on all three instances, one found and proved the optimum of 226 on \texttt{behnke31}. In an assisted session, the agent wrote a reusable tabu search, which found 406 on \texttt{behnke16} in 56.5\,s with 32 threads and the optimum on \texttt{behnke31}, but did not improve \texttt{behnke14} (232 against 230).

\begin{table}[!htbp]
\caption{All agent runs on the Behnke instances}
\label{tab:behnke}
\footnotesize
\begin{tabular*}{\textwidth}{@{\extracolsep{\fill}}l r r r@{}}
\toprule
 & \texttt{behnke14} & \texttt{behnke16} & \texttt{behnke31} \\
\midrule
FJSPLib LB / UB & 221 / 230 & 391 / 412 & 226 / 226 \\
\midrule
Auto, 5 sessions (A) & -- & \textbf{404}, \textbf{411}, 415, 417, 419 & -- \\
Auto, 5 more (A) & -- & 423, 430, 435, 450, 479 & -- \\
Auto, session 6 (A) & 235 & 416 & 229 \\
Auto, session 7 (A) & 234 & 415 & 226 \\
Assisted solver (B) & 232 & \textbf{406} & 226 \\
\bottomrule
\end{tabular*}
\end{table}

\paragraph{Generalization to the full Behnke set} \citet{behnke2012} published 60 instances, of which FJSPLib lists 32 as open and 28 as closed at the commit that we have used (f89e5b0, checked on October 11, 2026). We have run the reusable tabu search of the agent on all 60 without any change. As baselines, we have also run standard models in CP-SAT and in CP Optimizer~22.2 with the same budget, with one run of 120\,s with 8 threads per instance and method, on the same workstation. Both models use one optional interval for each eligible machine of an operation, which is the usual model for the flexible job shop. The tabu search and CP-SAT runs ran three at a time on the 16 cores, and the CP Optimizer runs two at a time, so CP Optimizer had slightly more cores per run. In a second block, we have given the tabu search one run of 1{,}200\,s with 16 threads on each open instance. These runs ran two at a time, so each had about 8 cores, half the cores of the original \texttt{behnke16} run (32 threads on the 16 cores). CP Optimizer received the same 1{,}200\,s with 16 threads, but one run at a time, so it had twice the cores of the tabu search in this block. Hexaly~15.0 ran with its standard flexible job-shop model, the same time limits and threads, and one run at a time. An independent checker verified all 336 schedules. Each value comes from a single run per instance and method. Differences of one or two units on a single instance are within the run-to-run variation (the 1{,}200\,s run of the tabu search ended worse than its 120\,s run on \texttt{behnke55} and \texttt{behnke59}). Table~\ref{tab:behnkegen} shows the results, and the public repository contains the values of all 60 instances. The table lists the open instances on which the tabu search or CP-SAT has reached or improved the best-known upper bound of FJSPLib \citep{fjsplib}. With 1{,}200\,s, CP Optimizer also reached it on \texttt{behnke15} and \texttt{behnke52}. TTB is the time in seconds to the best value of the 1{,}200\,s run of the tabu search, bold values are below the best-known value, and the superscript $\mathrm{o}$ marks a value that CP-SAT has reported as optimal.

\begin{table}[!htbp]
\caption{Open Behnke instances reached or improved by the tabu search or CP-SAT}
\label{tab:behnkegen}
\footnotesize
\setlength{\tabcolsep}{3pt}
\begin{tabular*}{\textwidth}{@{\extracolsep{\fill}}l r r r r r r r r r r@{}}
\toprule
 & & & \multicolumn{2}{c}{Agent tabu} & & & \multicolumn{2}{c}{CP Optimizer 22.2} & \multicolumn{2}{c}{Hexaly 15.0} \\
\cmidrule(lr){4-5}\cmidrule(lr){8-9}\cmidrule(lr){10-11}
Instance & LB & UB & 120\,s & 1{,}200\,s & TTB & CP-SAT & 120\,s & 1{,}200\,s & 120\,s & 1{,}200\,s \\
\midrule
\texttt{behnke12} & 213 & 219 & 222 & 219 & 172 & 243 & 225 & 221 & 233 & 227 \\
\texttt{behnke13} & 223 & 229 & 230 & \textbf{227} & 166 & 254 & 238 & 235 & 246 & 239 \\
\texttt{behnke16} & 391 & 412 & \textbf{409} & \textbf{406} & 641 & 497 & 423 & 416 & 432 & 422 \\
\texttt{behnke32} & 220 & 224 & 226 & \textbf{223} & 28 & 238 & 225 & \textbf{223} & 230 & 229 \\
\texttt{behnke34} & 219 & 223 & 224 & 223 & 4 & 245 & \textbf{222} & \textbf{221} & 229 & 225 \\
\texttt{behnke35} & 211 & 214 & 216 & \textbf{213} & 78 & 222 & 215 & \textbf{213} & 220 & 219 \\
\texttt{behnke40} & 415 & 419 & \textbf{418} & \textbf{417} & 307 & 495 & 425 & \textbf{418} & 432 & 420 \\
\texttt{behnke48} & 124 & 125 & 125 & 125 & 0 & 125$^{\mathrm{o}}$ & 125 & 125 & 126 & 125 \\
\texttt{behnke53} & 211 & 215 & \textbf{214} & \textbf{213} & 46 & \textbf{212}$^{\mathrm{o}}$ & 216 & \textbf{214} & 219 & 217 \\
\texttt{behnke55} & 221 & 223 & \textbf{222} & 223 & 336 & 228 & 224 & \textbf{222} & 226 & 226 \\
\texttt{behnke56} & 384 & 390 & 392 & \textbf{388} & 681 & 429 & 397 & 392 & 401 & 393 \\
\texttt{behnke57} & 385 & 390 & 391 & \textbf{389} & 41 & 450 & 399 & 396 & 402 & 392 \\
\texttt{behnke59} & 392 & 398 & \textbf{394} & \textbf{395} & 1195 & 458 & 400 & 398 & 405 & 401 \\
\texttt{behnke60} & 397 & 402 & 402 & \textbf{400} & 115 & 446 & \textbf{401} & \textbf{401} & 406 & 405 \\
\bottomrule
\end{tabular*}
\end{table}

With 120\,s, the tabu search went below the best-known upper bound on 5 of the 32 open instances, matched it on 2 more, and found the optimum on 18 of the 28 closed instances. On average, it was 0.60\% above the best-known upper bound, against 1.42\% for CP Optimizer, 2.85\% for Hexaly, and 6.42\% for CP-SAT. It was better than CP Optimizer on 31 instances, equal on 21, and worse on 8. With 1{,}200\,s on the open instances, the tabu search went below the best-known upper bound on 10 of them and matched it on 4 more, and it was 0.24\% above it on average, against 0.79\% for CP Optimizer with all 16 cores and 2.01\% for Hexaly. It was better than CP Optimizer on 19 open instances, equal on 9, and worse on 4. CP Optimizer went below the best-known upper bound on 7 open instances and found a new best-known value of 221 on \texttt{behnke34}. Hexaly improved no best-known value. A re-run of the tabu search with the same cores as CP Optimizer gave 0.45\% and 0.34\%. Over both blocks, the tabu search improved the best-known upper bound of 11 of the 32 open instances. CP-SAT found the optimum on 19 closed instances and reported 212 on \texttt{behnke53} and 125 on \texttt{behnke48} as optimal. Following the rule of Section~\ref{sec:verif}, we do not count these verdicts. However, FJSPLib now lists \texttt{behnke48} as closed and \texttt{behnke53} with a lower bound of 212, so our CP-SAT schedule of 212 closes \texttt{behnke53} with respect to the current library. On the large instances, CP-SAT was often 10--23\% above the best-known value.

\paragraph{Improvement campaign (mode~C)} Starting from the two Behnke solvers of the agents and the schedules above, one Claude Fable~5.1 session worked for 13.5\,h on our workstation, with open-source tools only and the goal of closing as many instances as possible. It merged the two solvers into a new local search and used CP-SAT to search for schedules with a given makespan. The campaign reached the FJSPLib lower bound on six open instances (\texttt{behnke32}, 34, 40, 55, 56, and 60), which are therefore closed, and on all 28 closed instances. Thus, 34 of the 60 instances now have a schedule at the lower bound. Its values are below the best-known upper bound on 18 of the 32 open instances, against 11 for modes A and B. The new local search produced most of the improvements on the 100-job instances, and eight of the new values, including the six closures, came from CP-SAT runs that the agent started with a fixed or capped makespan (Table~\ref{tab:modec}). Other processes ran on the workstation at the same time, so we do not compare run times. FJSPLib has changed during this study. Its current version (commit 8a6e36e, 22 October 2026) lists 27 open instances. Another contributor reported the same values as this campaign for \texttt{behnke52} and \texttt{behnke60} on 12 October, one week before the campaign, and CP-SAT runs raised several lower bounds. With respect to this version, the values of this campaign remain below the upper bound on 16 instances and close six (\texttt{behnke32}, 34, 38, 40, 55, and 56). On 9 November 2026, a further report to FJSPLib has claimed values equal to ours on \texttt{behnke32}, 34, 40, and 55 and better on \texttt{behnke35}, 53, 57, and 59 (on \texttt{behnke53}, equal to our CP-SAT value of 212), whereas our 387 on \texttt{behnke38} is better than its 388. 

\subsection{Flexible job shop with transport robots}\label{sec:res-fjspt}

\begin{figure}[!htbp]
  \centering
  \includegraphics[width=\textwidth]{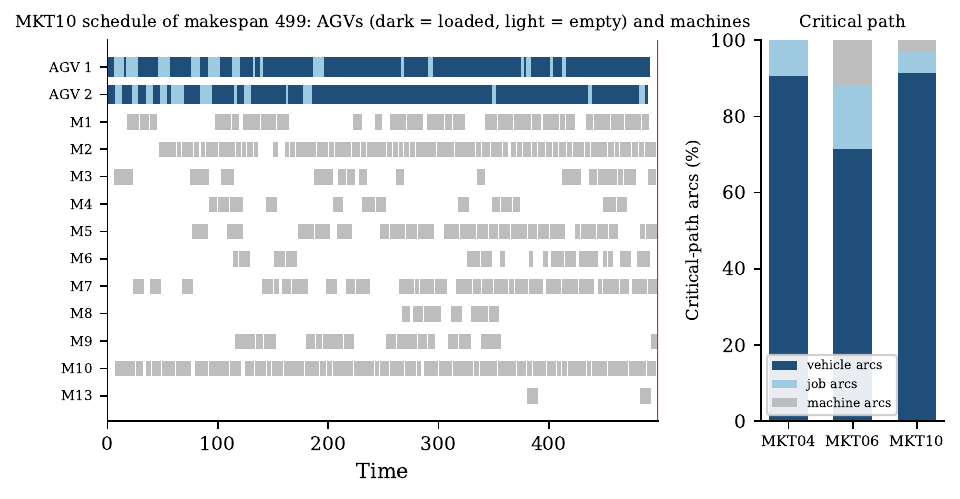}
  \caption{MKT10 schedule with makespan 499 and critical-path arcs}
  \label{fig:fjsptanatomy}
\end{figure}

\paragraph{Lower bounds} We have re-implemented the bounds of Section~\ref{sec:fjsptlb} independently. CP-SAT and CPLEX proved the same optimal values of bound~(1) on all 51 instances, and CP-SAT and CP Optimizer those of relaxation~(3) on all Chambers--Barnes, small, and medium instances, before any rounding. No bound is above any of the 177 schedules certified by our checker, the 40 replayed schedules of \citet{han2024dcga} and \citet{meng2025cp}, or the 18 proven optima, and no bound exceeds the optimum of 403 random and sub-sampled instances. The bounds of \citet{homayouni2021}, recomputed on our files, reproduce their published values except on MKT03 and MKT04, which suggests that their experiments on these two instances used other data. An earlier relaxation written by an agent in this study is a special case of bound~(1). Details are in the public repository.
\paragraph{Structure of the solutions} In the best MKT10 schedule, with a makespan of 499, the two vehicles are never idle. They carry jobs for 809 time units and drive empty for 172.5 units, so the makespan is half of this work plus a tail of about 8 units. On the critical path, 64 of 71 arcs connect consecutive tasks of the same vehicle when ties are broken toward job arcs (Figure~\ref{fig:fjsptanatomy}). The left panel of this figure shows the two vehicles (dark for loaded travel, light for empty travel) and the machines over time. The right panel shows the share of vehicle, job, and machine arcs on the critical path of the best schedules of MKT04, MKT06, and MKT10.

\paragraph{Progress over time} During the assisted session, the first runs on MKT10, 20 minutes into the session, found 541.5, already below all values published before 2025. The 2-opt$^\ast$ route exchange gave the largest gain of any single move (548 to 537.5 at equal budget). Slower cooling then yielded 522, ruin-and-recreate 518, and long independent runs 499 after 2.9 hours. The run that found 499 used one thread for 31 minutes. By then, the campaign had used about 40 thread-hours on MKT10.

\begin{table}[!htbp]
\caption{Results on the ten MKT instances with two AGVs}
\label{tab:mkt}
\scriptsize
\setlength{\tabcolsep}{1.6pt}
\begin{tabular*}{\textwidth}{@{\extracolsep{\fill}}l l r r r r r r r r r r r r r@{}}
\toprule
 & & & \multicolumn{6}{c}{Published upper bounds} & & & \multicolumn{3}{c}{Protocol (10 runs)} & \\
\cmidrule(lr){4-9}\cmidrule(lr){12-14}
Instance & Size (ops) & LB & LAHC & IGA & TS & DRL & GA & CP & Ours & Change & Range & Median & Below & CP-SAT \\
\midrule
MKT01 & $10\times6$ (55) & 167 & 187 & 177 & 187 & 153$^{\times}$ & 175 & \underline{174} & \textbf{173} & $-0.6$ & 173--175 & 174 & 1/10 & 184 \\
MKT02 & $10\times6$ (58) & 120 & 148 & 126 & 134 & \underline{122} & 124 & 124 & 124 & $+1.6$ & 124 & 124 & 0/10 & 156 \\
MKT03 & $15\times8$ (150) & 258 & 371 & 342 & 342 & 356 & 322 & \underline{313} & \textbf{294} & $-6.1$ & 294--298 & 297 & 10/10 & 556 \\
MKT04 & $15\times8$ (90) & 235 & 225$^{\times}$ & 244 & 254 & 213$^{\times}$ & 244 & \underline{237} & 238 & $+0.4$ & 238--240 & 239.5 & 0/10 & 310 \\
MKT05 & $15\times4$ (106) & 285 & 312 & 295 & 310 & 251$^{\times}$ & 290$^{\ddagger}$ & 286$^{\ddagger}$ & \textbf{286} & $-2.4$ & 286--288 & 288 & 10/10 & 397 \\
MKT06 & $10\times15^{\mathrm{c}}$ (150) & 244.5 & 389.5 & 321.5 & 336.5 & 329 & 300 & \underline{288.5} & \textbf{279.5} & $-3.1$ & 279.5--284 & 281.5 & 10/10 & 457 \\
MKT07 & $20\times5$ (100) & 228 & 291 & 267 & 289 & \underline{243} & 260 & 248 & 245 & $+0.8$ & 245--247 & 246 & 0/10 & 363 \\
MKT08 & $20\times10$ (225) & 662.5 & 846 & 780.5 & 770.5 & 723 & 757.5 & \underline{722} & \textbf{696} & $-3.6$ & 696--709 & 702.25 & 10/10 & 1184.5 \\
MKT09 & $20\times10$ (240) & 566.5 & 794 & 715 & 709 & 681 & 663 & \underline{632} & \textbf{598} & $-5.4$ & 601--610.5 & 606 & 10/10 & 1233 \\
MKT10 & $20\times15^{\mathrm{c}}$ (240) & 456 & 712.5 & 645 & 616.5 & 657 & 549 & \underline{525} & \textbf{499} & $-5.0$ & 505--515 & 509.5 & 10/10 & 1080 \\
\bottomrule
\end{tabular*}
\tablenote{\tabnotemark{c}Header value of the instance file. MKT06 uses machines 1--10, and MKT10 uses machines 1--10 and 13. Bold: below the reference value. $^{\times}$Below LB. $^{\ddagger}$The published schedule assigns an operation to a machine on which it is not eligible.}
\end{table}

\paragraph{MKT instances} Table~\ref{tab:mkt} gives the results on the ten MKT instances. LB is bound~(1) of Section~\ref{sec:fjsptlb}, which is above the best published lower bound of \citet{homayouni2021} on all ten instances. The published upper bounds are those of LAHC \citep{homayouni2021}, IGA \citep{meng2023}, tabu search (TS) \citep{berterottiere2024}, HA-DQN (DRL) \citep{dong2025hadrl}, the better of DCGA \citep{han2024dcga} and IVNDGA \citep{jia2025ivndga} (GA), and the better of the CP Optimizer models and the CP-assisted genetic algorithm of \citet{meng2025cp} (CP). The ant colony hybrid of \citet{he2025aco} is never the smallest. A published value below LB cannot be the makespan of a feasible schedule on these instance files, so we do not use it as a reference. We replayed every published schedule with our checker. All reproduce their makespans and are feasible, except on MKT05, where the schedules of 290 \citep{han2024dcga} and 286 \citep{meng2025cp} assign an operation to a machine on which it is not eligible in the distributed file. The next value there is 293, from the CP Optimizer models of \citet{meng2025cp}. If these two schedules were accepted, MKT05 would be a tie, and six of the ten instances would be improved.  The reference value is the smallest remaining published value. It is underlined, except on MKT05, where it is 293 and not shown. ``Change'' is relative to it, and ``Below'' counts the runs of the protocol (ten independent 20-minute single-thread runs per instance, 28 concurrent processes on 32 threads) that are below it. The last column is one 10-minute CP-SAT run with 8 threads from a greedy start. The lower bounds that CP-SAT proved in these runs (61 to 587) are below LB on all ten instances.

On seven of the ten instances (MKT01, MKT03, MKT05, MKT06, MKT08, MKT09, and MKT10), the certified makespans are below the reference value. The largest changes are $-6.1\%$ on MKT03, $-5.4\%$ on MKT09, and $-5.0\%$ on MKT10. On all ten instances, the agent solver is at least as good as the CP Optimizer models of \citet{meng2025cp}, which ran for 3{,}600\,s, and it is better on nine. The CP-assisted genetic algorithm of the same study reports 237 on MKT04, one unit below our 238, without a published schedule for this instance. The HA-DQN values on MKT01, MKT04, and MKT05 (153, 213, and 251) and the LAHC value on MKT04 (225) are below LB, and the closed-form bound of Section~\ref{sec:fjsptlb} (165, 234, and 282) is sufficient to show this. If one counts all published values without these checks, our values are below all of them on five instances. On MKT02 and MKT07, \citet{dong2025hadrl} report 122 and 243, below our 124 and 245. These two values are not below LB, and they were published without schedules, so we could not check them. With respect to LB, our best schedules are 0.4--14.3\% above it. On MKT05, the schedule of 286 is only one unit above the bound of 285. As a fixed protocol, we have run the final solver ten times for 20 minutes on each instance (Table~\ref{tab:mkt}). On the seven instances just mentioned, 61 of the 70 runs are below the reference value.

\paragraph{Chambers--Barnes instances} Table~\ref{tab:cb} gives the results on the 21 Chambers--Barnes transport instances. In this table, LB is relaxation~(3) of Section~\ref{sec:fjsptlb}, which CP-SAT and CP Optimizer solve to optimality on all 21 instances. ``L/I/D'' is the smallest value of LAHC \citep{homayouni2021}, IGA \citep{meng2023}, and DCGA \citep{han2024dcga}, ``DRL'' the value of HA-DQN \citep{dong2025hadrl}, and ``CP'' the best value of \citet{meng2025cp}. A dagger marks a CP value without a published schedule. All other CP values are the makespans of the published schedules replayed with our checker, which differ from the reported values on \texttt{mt10xyzt} (901 instead of 916), \texttt{seti5xxxt} (1{,}316 instead of 1{,}318), and \texttt{seti5cct}, whose schedule omits a transport of zero length and needs 1{,}263 instead of 1{,}261 when it is inserted. The column $\Delta$ is relative to the smallest published value that is not below LB. The agent-written solver is better than L/I/D on 12 of the 21 instances. However, the approach of \citet{meng2025cp}, based on CP Optimizer, is better than our solver on 20 instances and equal on \texttt{mt10xyt}, where both schedules are optimal. Our solver is 0.0--6.5\% above the best published value. It is weaker on these instances because processing times are much longer than travel times, and the critical path consists of machine arcs only, which a routing neighborhood cannot shorten. The relaxation remains informative. Seven HA-DQN values, published without schedules, are below LB, among them 872 and 816 on \texttt{mt10xxxt} and \texttt{mt10xyzt}, which are also below the optima of the same instances without transport in FJSPLib \citep{fjsplib}. \citet{meng2025cp} report 14 of the 21 instances as proven optimal. On seven of them, their published schedule reaches LB (asterisks in Table~\ref{tab:cb}), which confirms optimality independently of their solver. On \texttt{mt10xyzt}, under the benchmark model, which the schedules of \citet{meng2023} and \citet{berterottiere2024} reproduce exactly, their schedule has a makespan of 901, equal to LB. The reported 916 may come from a different evaluation.

\begin{table}[!htbp]
\caption{Chambers--Barnes instances with two AGVs}
\label{tab:cb}
\footnotesize
\setlength{\tabcolsep}{3pt}
\scriptsize
\setlength{\tabcolsep}{2pt}
\begin{tabular*}{\textwidth}{@{\extracolsep{\fill}}l r r r r r r l r r r r r r@{}}
\toprule
Instance & LB & L/I/D & DRL & CP & Ours & $\Delta$\% & Instance & LB & L/I/D & DRL & CP & Ours & $\Delta$\% \\
\midrule
\texttt{mt10c1t} & 978 & 991 & 969$^{\times}$ & 979 & 996 & $+1.7$ & \texttt{setb4xxxt} & 961 & 992 & -- & 961$^{\ast}$ & 978 & $+1.8$ \\
\texttt{mt10cct} & 948 & 973 & 933$^{\times}$ & 948$^{\ast}$ & 973 & $+2.6$ & \texttt{setb4xyt} & 943 & 963 & -- & 948 & 975 & $+2.8$ \\
\texttt{mt10xt} & 975 & 986 & 942$^{\times}$ & 978 & 987 & $+0.9$ & \texttt{setb4xyzt} & 939 & 972 & 932$^{\times}$ & 939$^{\ast}$ & 968 & $+3.1$ \\
\texttt{mt10xxt} & 979 & 991 & -- & 979$^{\ast}$ & 992 & $+1.3$ & \texttt{seti5c12t} & 1240 & 1347 & 1217$^{\times}$ & 1283$^{\dagger}$ & 1307 & $+1.9$ \\
\texttt{mt10xxxt} & 963 & 981 & 872$^{\times}$ & 963$^{\ast}$ & 990 & $+2.8$ & \texttt{seti5cct} & 1194 & 1334 & 1209 & 1263 & 1287 & $+6.5$ \\
\texttt{mt10xyt} & 950 & 959 & -- & 950$^{\ast}$ & 950$^{\ast}$ & $+0.0$ & \texttt{seti5xt} & 1245 & 1355 & -- & 1312$^{\dagger}$ & 1332 & $+1.5$ \\
\texttt{mt10xyzt} & 901 & 928 & 816$^{\times}$ & 901$^{\ast}$ & 922 & $+2.3$ & \texttt{seti5xxt} & 1253 & 1373 & -- & 1322$^{\dagger}$ & 1326 & $+0.3$ \\
\texttt{setb4c9t} & 954 & 1004 & 958 & 957 & 1005 & $+5.0$ & \texttt{seti5xxxt} & 1246 & 1350 & -- & 1316 & 1330 & $+1.1$ \\
\texttt{setb4cct} & 940 & 979 & 946 & 943 & 991 & $+5.1$ & \texttt{seti5xyt} & 1195 & 1346 & -- & 1258$^{\dagger}$ & 1285 & $+2.1$ \\
\texttt{setb4xt} & 959 & 994 & -- & 961 & 1002 & $+4.3$ & \texttt{seti5xyzt} & 1190 & 1305 & 1208 & 1240 & 1274 & $+5.5$ \\
\texttt{setb4xxt} & 957 & 1009 & -- & 961 & 993 & $+3.3$ &  & & & & & &  \\
\bottomrule
\end{tabular*}
\tablenote{$^{\times}$Below LB. $^{\dagger}$No published schedule. $^{\ast}$Equal to LB, hence optimal.}
\end{table}

\paragraph{Generalization to other FJSPT instances} We have run the final transport solver without any change on the other two-vehicle sets of \citet{homayouni2021}: ten small instances (SFJST) and ten medium instances (MFJST, built from the instances of \citealt{fattahi2007}). Before the runs, we have checked the instance files, and the published schedules of \citet{berterottiere2024} and \citet{meng2023} replay to their published makespans on these files. Each instance received five runs of 60\,s on one thread. On the small instances, 48 of the 50 runs found the proven optimum, and every instance was solved to optimality by at least four runs. The medium instances are less uniform. The solver found the proven optimum \citep{homayouni2021} on MFJST01, 03, and 06 in at least one of five runs, and matched the best published value of 900 on MFJST08, which equals relaxation~(3) and is therefore optimal. On MFJST02, 04, 05, and 07, its best run stayed 0.5--2.3\% above the optimal value. On MFJST09, it found 1{,}090, below the values of DCGA and IVNDGA \citep{han2024dcga,jia2025ivndga} but above the 1{,}075 of \citet{meng2025cp}, whose schedule we replayed and which equals relaxation~(3), so it is optimal. On MFJST10, it stayed 1.4\% above the best published value of 1{,}227. These instances have long processing times compared with the travel times, and the final temperatures of the annealing (0.2 to 4 time units) are too small for them. We have not tuned the solver for these instances, because the goal was to test the unchanged method.

\paragraph{Open-source baseline} As an open-source baseline, we have run CP-SAT on a complete constraint model of the problem. The model uses optional intervals for the machine choice and a multiple-circuit constraint for the routes of the two vehicles, with empty travel times that depend on the chosen machines. It accepts every published schedule and every schedule of this work as feasible. CP-SAT started from a greedy dispatch schedule and ran for 10 minutes with 8 threads on the same workstation. Table~\ref{tab:mkt} compares its result with the ten 20-minute single-thread runs of the agent solver. Although the agent solver used a quarter of the CPU time (1{,}200 thread-seconds against 4{,}800), CP-SAT ended 6--112\% above the median of the agent solver.

\subsection{Summary of the comparison with solvers}\label{sec:res-compare}

Table~\ref{tab:compare} collects the comparisons of Sections~\ref{sec:res-jssp}--\ref{sec:res-fjspt}. Gaps are relative to the best-known lower bound, and the Behnke rows show the mean deviation from the best-known upper bound. The upper part contains only runs on our workstation with the same budget for both solvers. The lower part compares with published values from other hardware. In the Behnke rows with CP Optimizer, values in parentheses come from the re-run of the tabu search with the same cores.

\begin{table}[!htbp]
\caption{Summary of the comparisons with other solvers}
\label{tab:compare}
\scriptsize
\setlength{\tabcolsep}{2pt}
\begin{tabular*}{\textwidth}{@{\extracolsep{\fill}}l l l l l l@{}}
\toprule
 & \multicolumn{3}{c}{Reference solver} & \multicolumn{2}{c}{Agent-written solver} \\
\cmidrule(lr){2-4}\cmidrule(lr){5-6}
Instances & Solver & Budget & Result & Budget & Result \\
\midrule
\multicolumn{6}{@{}l}{\emph{Same workstation, same budget}} \\
\texttt{tai\_100\_1000\_1..10} & CP Optimizer 22.2 & 10 min, 1 core & gap 4.0--6.6\% & 10 min, 1 core & gap 0.0--1.8\% \\
\texttt{tai\_100\_1000\_1..10} & CP Optimizer 22.2 & 30 min, 1 core & gap 3.7--5.8\% & 30 min, 1 core & gap 0.0--1.8\% \\
\texttt{tai\_100\_1000\_1..10} & CP Optimizer 22.2 & 30 min, 16 cores & gap 2.6--5.3\% & 30 min, 16 cores & gap 0.0--1.6\% \\
\texttt{tai\_100\_1000\_1..10} & Hexaly 15.0 & 10 min, 1 core & gap 2.7--4.5\% & 10 min, 1 core & gap 0.0--1.8\% \\
\texttt{tai\_100\_1000\_1..10} & Hexaly 15.0 & 30 min, 1 core & gap 2.7--4.4\% & 30 min, 1 core & gap 0.0--1.8\% \\
\texttt{tai\_100\_1000\_1..10} & Hexaly 15.0 & 30 min, 16 cores & gap 2.4--4.2\% & 30 min, 16 cores & gap 0.0--1.6\% \\
Behnke 1--60 & CP Optimizer 22.2 & 120\,s, 8 threads & mean +1.42\% & 120\,s, 8 threads & mean +0.60\% (+0.45\%) \\
Behnke, 32 open & CP Optimizer 22.2 & 1{,}200\,s, 16 threads & mean +0.79\% & 1{,}200\,s, 16 threads & mean +0.24\% (+0.34\%) \\
Behnke 1--60 & Hexaly 15.0 & 120\,s, 8 threads & mean +2.85\% & 120\,s, 8 threads & mean +0.60\% \\
Behnke, 32 open & Hexaly 15.0 & 1{,}200\,s, 16 threads & mean +2.01\% & 1{,}200\,s, 16 threads & mean +0.24\% \\
Behnke 1--60 & CP-SAT (open source) & 120\,s, 8 threads & mean +6.42\% & 120\,s, 8 threads & mean +0.60\% \\
MKT01--10 & CP-SAT (open source) & 10 min, 8 threads & e.g.\ MKT10: 1{,}080 & 20 min, 1 thread & e.g.\ MKT10: 509.5 (median) \\
\midrule
\multicolumn{6}{@{}l}{\emph{Published values, other hardware}} \\
\texttt{tai\_100\_1000\_1..10} & OptalCP & 10 min, 4-core laptop & gap 2.6--3.8\% & 10 min, 1 core & gap 0.0--1.8\% \\
\texttt{tai\_100\_1000\_1} & OptalCP & 180\,min, 16 workers & 533{,}080 & 2.7--16.2\,s, 1 core & 533{,}000 \\
MKT01--10 & CP Optimizer models \citep{meng2025cp} & 3{,}600\,s & published values & Table~\ref{tab:mkt} & better on 9, equal on 1 \\
Chambers--Barnes & CP-based approach \citep{meng2025cp} & as published & published values & 150--1{,}200\,s, 1 thread & worse on 20, equal on 1 \\
\bottomrule
\end{tabular*}
\end{table}

\subsection{New best-known values and closed instances}\label{sec:res-newbest}

Table~\ref{tab:newbest} summarizes the instances on which the agents found a schedule better than the best-known value before this work (``New'') or equal to the lower bound (``Closed''), with the range of the time to solution of modes A and B. Appendix~\ref{app:newbest} lists all 28 instances of modes A and B with the previous bounds, the value, the run that found it, and its time to solution.

\begin{table}[!htbp]
\caption{New best-known values and closed instances by problem family}
\label{tab:newbest}
\footnotesize
\begin{tabular*}{\textwidth}{@{\extracolsep{\fill}}l r r r r r@{}}
\toprule
 & \multicolumn{3}{c}{Modes A and B} & \multicolumn{2}{c}{With mode C} \\
\cmidrule(lr){2-4}\cmidrule(lr){5-6}
Family & New & Closed & Time & New & Closed \\
\midrule
Job shop & 10 & 2 & 7 min--30 min & 10 & 4 \\
Behnke--Geiger & 11 & 0 & 28\,s--3.2 h & 18 & 6 \\
FJSPT MKT & 7 & 0 & 3 min--31 min & 7 & 0 \\
Chambers--Barnes & 0 & 0 & -- & 0 & 0 \\
FJSPT MFJST & 0 & 0 & -- & 0 & 0 \\
\midrule
Total & 28 & 2 & & 35 & 10 \\
\bottomrule
\end{tabular*}
\end{table}

Two instances are now closed. On \texttt{tai\_100\_1000\_2} and \texttt{\_3}, the corridor solver of the agent found a schedule whose makespan equals the lower bound given by the longest job. Of the 28 values, 27 come from the agent-written solvers in mode~B and one from an autonomous session (\texttt{behnke16}, 404). The agent-written solvers found their values in 28\,s to 31 minutes, whereas the autonomous session that found 404 needed about 3.2 hours. Values of our own CP-SAT and CP Optimizer baselines are not counted. Including the improvement campaigns of mode~C, the agents found new best-known values on 35 instances and closed ten (Tables~\ref{tab:newbest} and~\ref{tab:modec}). These counts refer to the bounds when our experiments started (FJSPLib commit f89e5b0, 11 October 2026). With respect to the current FJSPLib, the Behnke values of mode~C are new on 16 instances and close six (Section~\ref{sec:res-fjsp}). If the report of 9 November 2026 is accepted, they are new on 8 instances and close two. The new schedules will be submitted to JSPLib and FJSPLib.

\section{Discussion}\label{sec:disc}

\subsection{The agent as the solver, as the author of a solver, or as its improver}

Table~\ref{tab:modes} compares modes A and B of Section~\ref{sec:modes}.

In mode~A, the agent itself is the algorithm. The user needs neither optimization knowledge nor code, and the agent adapts to each instance, up to changing the formulation. The job-shop session moved to a fixed-target formulation with time windows, a known technique in exact job-shop methods \citep{carlier1994heads}, and applied it at a scale of 100{,}000 operations. However, every new schedule requires a new session, which on \texttt{tai\_100\_1000\_1} cost US\$2.43--22.69 and took 18--300 minutes, and the results vary considerably (527{,}000 to 554{,}332).

\begin{table}[!htbp]
\caption{Modes A and B of a coding agent}
\label{tab:modes}
\footnotesize
\renewcommand{\arraystretch}{1.2}
\begin{tabular*}{\textwidth}{@{\extracolsep{\fill}}>{\raggedright\arraybackslash}p{0.26\textwidth}>{\raggedright\arraybackslash}p{0.33\textwidth}>{\raggedright\arraybackslash}p{0.33\textwidth}@{}}
\toprule
 & \textbf{Mode A}: the agent is the solver & \textbf{Mode B}: the agent writes a solver \\
\midrule
\textit{Output} & one schedule & a reusable program \\
\textit{Development cost (once)} & none & one session: US\$10 (autonomous, 101.6\,min) or an assisted session (5--12\,h) \\
\textit{Cost of a new schedule} & a new session: US\$2--23 & CPU time only \\
\textit{Time to a new schedule} & 0.3--5\,h & seconds to 30 minutes \\
\textit{Run-to-run variance} & high & low; controlled by seeds \\
\textit{Repeatable} & no & yes \\
\textit{Similar instances} & new session each & runs unchanged \\
\textit{New structure} & adapts by itself & needs a new session \\
\textit{Auditable} & partly (transcript) & fully (source code) \\
\bottomrule
\end{tabular*}
\end{table}

In mode~B, the agent writes a solver once, and from then on a new schedule costs only CPU time. Because the solver is code, an expert can read it, test it, and change it. It is also repeatable. We have recompiled the transport solver from its source, have used the same seed, and have obtained a schedule identical to the MKT10 schedule of this paper. The solvers also apply to similar instances without change (Section~\ref{sec:results}). On the job-shop family, the solver from one session has improved all ten instances and has closed two of them. The reusable tabu search has improved the best-known upper bound of 11 open Behnke instances, and the transport solver has found every optimum on the small transport instances. The drawback is that a solver tailored to one structure can be weak on another. The transport solver is strong when the vehicles are the bottleneck and weaker when the machines are (Section~\ref{sec:res-fjspt}).

In practice, the two modes work best together. An agent session identifies the structure of a problem and produces a method. A short assisted session then turns this method into a solver that the user runs again whenever a new instance with a similar structure arrives. This matches a common situation in industry, where a plant solves the same type of problem every day with different data.

Mode~C follows from these results. Once the solvers of mode~B had been applied to whole instance families without change, we studied whether an agent can take the next step without human input: start from everything that modes A and B produced and coordinate its own campaign. In both campaigns, the agent planned the experiments, changed or combined the solvers, distributed the runs over the cores, and abandoned some search directions. We only set the goal and stopped the session. On the job shop, it analyzed why the corridor solver stalls close to the lower bound and changed the search instead of adding compute. On the Behnke instances, it merged the two agent-written solvers and used CP-SAT at a fixed makespan. The three modes thus form a sequence: an agent produces a method, the method becomes a reusable solver, and an agent then improves this solver with the results obtained so far. A campaign needs more time and compute than a single run (24 and 13.5 hours), but still no expert work and no commercial software.

\subsection{Quality, speed, and cost}

In terms of quality, Section~\ref{sec:res-compare} gives the numbers. The published values of the commercial solvers come from other hardware and from time limits between 10 minutes and 6 hours. Our own runs of CP Optimizer and Hexaly remove this difference. With the same budget on the same workstation, both were weaker than the agent solvers on the job shop (2.4--5.3\% against 0.0--1.6\% above the lower bound with 16 cores), and weaker on the Behnke instances, where CP Optimizer nevertheless found the new best-known value 221 on \texttt{behnke34}. On the transport instances, the agent solver beats a published approach based on CP Optimizer when the vehicles are the bottleneck (MKT), but not when the machines are (Chambers--Barnes). Some agent results used more total CPU time than the reference runs, and the solver budgets do not include the time needed to develop the agent solvers (see below).

In terms of speed, the agent-written solvers are also fast. On one core, the corridor solver reaches the 180-minute OptalCP value of \texttt{tai\_100\_1000\_1} in 2.7--16.2\,s, the reusable tabu found 406 on \texttt{behnke16} after 56.5\,s, and the first runs of the transport session were below all MKT10 values published before 2025 after 20 minutes. Developing a solver takes longer: 101.6 minutes for the corridor solver, and about 5 and 12 hours for the two recorded assisted sessions.

In terms of cost, we separate the cost of developing a solver from the cost of a new schedule. The 16 autonomous job-shop sessions took 48 hours in total and cost US\$167 in API charges for the 15 with a recorded cost (median US\$12.72). The session that wrote the corridor solver cost US\$10.07. The interactive sessions ran under a Claude Max subscription (US\$100 per month) on a few days, together less than one month of subscription time, so we estimate their cost at about US\$100. The two campaigns of mode~C also ran within such a plan. At API list prices, the logged token use of the two recorded interactive sessions would cost about US\$199. Our own time for goals, review, and verification is not included. Once a solver exists, each new schedule costs only CPU time on commodity hardware and needs no license, whereas a developer subscription of CP Optimizer starts at US\$320 per user and month \citep{ibmcplexpricing}. A general-purpose solver must work well on every instance, whereas the agent can specialize to one structure. This specialization appears to be the main source of its advantage.

\subsection{Limitations}

We have tested only models from one provider (Anthropic), namely Claude Opus and Claude Sonnet, and Claude Fable in one campaign of mode~C. Controlled comparisons across models and problem classes are needed to show which models are suitable for this task (Section~\ref{sec:conc}).

All instances are benchmarks from the literature. Their generators can give them regular structure, for example the long, almost idle machines of the $100\times1000$ job shops or the clear stages of \texttt{behnke16}, and the agents have exploited such structure. It should be noted, however, that our conclusions may not hold for industrial problems, which are often less regular and include additional constraints.

The agents did not raise the lower bound of any job-shop or Behnke instance. On the transport instances, the agent of the assisted session wrote a relaxation that bound~(1) of Section~\ref{sec:fjsptlb} generalizes. The closed instances (\texttt{tai\_100\_1000\_2}, \texttt{\_3}, \texttt{\_4}, and \texttt{\_9}, and six Behnke instances in mode~C) are closed because the schedules reach the known lower bound. \texttt{tai\_100\_1000\_1} and \texttt{behnke16} remain open.

The generalization runs used short budgets and few seeds (three runs of 30 minutes, or single runs of 120\,s and 1{,}200\,s). We ran CP Optimizer and Hexaly with standard models and default parameters, as a typical user would. Models or parameters tuned by an expert for these instances may yield better results. On the medium transport instances, the unchanged solver stayed above the best known value on six of ten instances.

MKT10 served as the main development instance of the transport solver, which may favor it. For this reason, we report the full instance sets. On the Chambers--Barnes instances, whose critical paths consist of machine arcs, the transport solver is weaker than a recent approach based on CP Optimizer \citep{meng2025cp}. A solver tailored to the vehicles does not transfer to instances in which the machines are the bottleneck. Finally, \citet{dong2025hadrl} report lower values than ours on five MKT instances. Three of these values are below our lower bound. They do not publish schedules, so we could not check the other two with our checker.

\section{Conclusions and future research}\label{sec:conc}

We have examined whether a frontier LLM coding agent can replace a commercial solver on difficult scheduling problems. On three established problem families, the agents found schedules better than the best published values, using only open-source tools and their own code. In autonomous sessions, without expert knowledge, they improved the 100{,}000-operation job shop \texttt{tai\_100\_1000\_1} and the flexible job shop \texttt{behnke16}. The solvers they wrote transfer to similar instances: unchanged, the job-shop solver improved all ten instances of its family and closed two, and in mode~C agents closed two further job shops and six flexible job shops. With transport robots, the agent solver improved seven of the ten MKT instances but none of the Chambers--Barnes instances, on which a CP Optimizer approach is better. New, validated lower bounds for this problem are above eleven published values on the distributed instance files and certify eight published schedules as optimal. On large benchmark instances with exploitable structure, solvers written by current LLM agents can match or beat commercial solvers, both their published values and our own runs of CP Optimizer and Hexaly with the same budget. The cost is low. The job-shop solver took one autonomous session of 101.6 minutes and US\$10, and each new schedule then costs only CPU time. The agents are strongest when an instance has a structure that a general-purpose solver does not exploit.

This study provides evidence for models from one provider, three model lines, and three problem families, which were selected because they are open, large, and structurally rich. It does not show how other models perform, or how the results transfer to other classes of optimization problems. More studies are therefore required that compare different models on a wide range of optimization problems under controlled conditions. Such studies should vary the model family and size, keep the task text, the tools, the compute resources, and the time budget fixed, and report all runs, so that differences in performance can be attributed to the models.

A natural next step is a standardized benchmark for this purpose. It would collect established instances from many problem classes, for example the job shop, the flexible job shop, the open shop \citep{gonzalez1976openshop}, scheduling with transport resources, nurse rostering \citep{ceschia2019inrc2}, and timetabling \citep{mccollum2010itc2007}, each with its source, its best-known bounds, and an independent checker. Runs of several model families on the same benchmark would show how each model performs on each problem class, and the evaluation of new model versions on the same tasks would show how this performance changes over time. A benchmark that covers many problem classes would also reduce the risk of an instance selection that favors one method, which is a limitation of the present study.

A second open question is which behavior of the agents yields good solutions. Here, a change of formulation and the use of instance structure were decisive, but on a small number of sessions. A systematic analysis of many sessions could test whether these strategies are general.

\section*{Data and materials availability}

All materials of this study are available in a public repository at \url{https://github.com/andrealbrechtl/agent-heuristics-jobshop-supplement/}. For each problem family, it contains the instance files with their sources and checksums, the certified schedules, an independent checker, the solver source code, the run commands, and the solver logs. It also contains the replays of published schedules \citep{meng2023,berterottiere2024,han2024dcga,meng2025cp}, the models, schedules, and logs of our CP Optimizer, Hexaly, and CP-SAT baselines, the campaigns of mode~C, and the scripts and results of the lower bounds of Section~\ref{sec:fjsptlb}. For \texttt{setb4xxx}, the distributed data file repeats the block of \texttt{seti5xxx}, and we used the original Chambers--Barnes data. The task texts of the autonomous job-shop and \texttt{behnke16} runs and the first message of the assisted transport-robot session are included verbatim. The other runs used the same or very similar task texts.  One script re-checks every value in this paper.

\section*{Use of generative AI}

The use of LLM coding agents to produce the experimental results is the subject of the study and is described in Section~\ref{sec:method}. During the preparation of the manuscript, the authors used Claude (Anthropic) to review the text and to rewrite parts of it. The lower bounds of Section~\ref{sec:fjsptlb} were derived and implemented by Claude Opus~5.5 (\texttt{claude-opus-5-5}, in Claude Code, November 2026) during the review of the manuscript. Bound~(1) generalizes a relaxation written earlier by the agent of the transport session. The authors re-implemented and validated them independently with two solvers. The authors have checked the proofs.  The authors reviewed and edited all content and take full responsibility for the publication.

\section*{Acknowledgments}

This work was partly supported by the Canadian Natural Sciences and Engineering Research Council (NSERC) under grant 2025-03964. This support is gratefully acknowledged. This study was financed in part by the Coordena\c{c}\~ao de Aperfei\c{c}oamento de Pessoal de N\'{\i}vel Superior -- Brasil (CAPES) -- Finance Code 001 and by the Conselho Nacional de Desenvolvimento Cient\'{\i}fico e Tecnol\'ogico (CNPq).

\section*{Conflict of interest}

The authors have no conflict of interest to declare.

\bibliographystyle{itor}
\bibliography{refs}

\appendix
\section{New best-known values found by the agents}\label{app:newbest}

Table~\ref{tab:newbestfull} lists every instance on which an agent of mode~A or~B found a new best-known value (Section~\ref{sec:res-newbest}), and Table~\ref{tab:modec} the values of mode~C below both the best value of modes A and B (A/B) and the previous best-known value. LB and UB are the best-known bounds before this work (JSPLib, FJSPLib, and, for MKT, the bound of \citet{homayouni2021} and the reference value of Table~\ref{tab:mkt}). ``Run'' gives the resources of the run that found the value (16 $\times$ 1 core: 16 single-core workers), and ``Time'' the wall time to the value (FJSPT: to the end of its cooling run). Notes: $a$ from file timestamps; $c$ one run of a campaign of about 40 thread-hours; $d$ our CP-SAT baseline later found 212 and reported it as optimal. ``Compaction'' starts every operation of the mode-B schedule as early as possible in its machine order, ``corridor search'' is the solver of Section~\ref{sec:res-jssp}, and ``new local search'' that of the Behnke campaign.

\begin{table}[!htbp]
\caption{All instances on which an agent found a new best-known value}
\label{tab:newbestfull}
\scriptsize
\setlength{\tabcolsep}{3pt}
\begin{tabular*}{\textwidth}{@{\extracolsep{\fill}}l r r r l l l r@{}}
\toprule
Instance & LB & UB & Ours & Status & Found by & Run & Time \\
\midrule
\multicolumn{8}{@{}l}{\emph{Job shop}} \\
\texttt{tai\_100\_1000\_1} & 522{,}298 & 533{,}080 & \textbf{525{,}515} & improved & Agent solver (B) & 16 $\times$ 1 core, 30 min & 29.5 min \\
\texttt{tai\_100\_1000\_2} & 530{,}375 & 538{,}067 & \textbf{530{,}375} & closed & Agent solver (B) & 1 core, 30 min & 22.6 min \\
\texttt{tai\_100\_1000\_3} & 530{,}560 & 538{,}757 & \textbf{530{,}560} & closed & Agent solver (B) & 1 core, 30 min & 7.0 min \\
\texttt{tai\_100\_1000\_4} & 527{,}101 & 534{,}746 & \textbf{528{,}156} & improved & Agent solver (B) & 16 $\times$ 1 core, 30 min & 26.6 min \\
\texttt{tai\_100\_1000\_5} & 517{,}728 & 529{,}580 & \textbf{525{,}865} & improved & Agent solver (B) & 16 $\times$ 1 core, 30 min & 20.3 min \\
\texttt{tai\_100\_1000\_6} & 522{,}907 & 534{,}969 & \textbf{524{,}812} & improved & Agent solver (B) & 16 $\times$ 1 core, 30 min & 29.2 min \\
\texttt{tai\_100\_1000\_7} & 522{,}537 & 534{,}974 & \textbf{527{,}000} & improved & Agent solver (B) & 1 core, 30 min & 8.8 min \\
\texttt{tai\_100\_1000\_8} & 526{,}428 & 535{,}757 & \textbf{528{,}473} & improved & Agent solver (B) & 16 $\times$ 1 core, 30 min & 29.9 min \\
\texttt{tai\_100\_1000\_9} & 528{,}097 & 536{,}993 & \textbf{529{,}470} & improved & Agent solver (B) & 16 $\times$ 1 core, 30 min & 28.5 min \\
\texttt{tai\_100\_1000\_10} & 521{,}766 & 529{,}918 & \textbf{526{,}000} & improved & Agent solver (B) & 16 $\times$ 1 core, 30 min & 16.1 min \\
\midrule
\multicolumn{8}{@{}l}{\emph{Behnke--Geiger}} \\
\texttt{behnke16} & 391 & 412 & \textbf{404} & improved & Auto session (A) & 4-CPU container & $\approx$3.2 h$^{a}$ \\
\texttt{behnke13} & 223 & 229 & \textbf{227} & improved & Agent solver (B) & 16 thr, 1{,}200\,s & 2.8 min \\
\texttt{behnke32} & 220 & 224 & \textbf{223} & improved & Agent solver (B) & 16 thr, 1{,}200\,s & 28\,s \\
\texttt{behnke35} & 211 & 214 & \textbf{213} & improved & Agent solver (B) & 16 thr, 1{,}200\,s & 78\,s \\
\texttt{behnke40} & 415 & 419 & \textbf{417} & improved & Agent solver (B) & 16 thr, 1{,}200\,s & 5.1 min \\
\texttt{behnke53} & 211 & 215 & \textbf{213}$^{d}$ & improved & Agent solver (B) & 16 thr, 1{,}200\,s & 46\,s \\
\texttt{behnke55} & 221 & 223 & \textbf{222} & improved & Agent solver (B) & 8 thr, 120\,s & 31\,s \\
\texttt{behnke56} & 384 & 390 & \textbf{388} & improved & Agent solver (B) & 16 thr, 1{,}200\,s & 11.4 min \\
\texttt{behnke57} & 385 & 390 & \textbf{389} & improved & Agent solver (B) & 16 thr, 1{,}200\,s & 41\,s \\
\texttt{behnke59} & 392 & 398 & \textbf{394} & improved & Agent solver (B) & 8 thr, 120\,s & 92\,s \\
\texttt{behnke60} & 397 & 402 & \textbf{400} & improved & Agent solver (B) & 16 thr, 1{,}200\,s & 115\,s \\
\midrule
\multicolumn{8}{@{}l}{\emph{FJSPT, MKT}} \\
\texttt{MKT01} & 125 & 174 & \textbf{173} & improved & Agent solver (B) & 1 thr, 1{,}200\,s & 3.4 min \\
\texttt{MKT03} & 206 & 313 & \textbf{294} & improved & Agent solver (B) & 1 thr, 1{,}200\,s & 12.2 min \\
\texttt{MKT05} & 231 & 293 & \textbf{286} & improved & Agent solver (B) & 1 thr, 900\,s & 13.3 min \\
\texttt{MKT06} & 222.8 & 288.5 & \textbf{279.5} & improved & Agent solver (B) & 1 thr, 1{,}200\,s & 13.2 min \\
\texttt{MKT08} & 591.5 & 722 & \textbf{696} & improved & Agent solver (B) & 1 thr, 1{,}200\,s & 10.9 min \\
\texttt{MKT09} & 481.8 & 632 & \textbf{598} & improved & Agent solver (B) & 1 thr, 2{,}700\,s & 21.2 min \\
\texttt{MKT10} & 386.8 & 525 & \textbf{499} & improved & Agent solver (B) & 1 thr, 1{,}900\,s & 30.6 min$^{c}$ \\
\bottomrule
\end{tabular*}
\end{table}

\begin{table}[!htbp]
\caption{Values of the two improvement campaigns of mode~C}
\label{tab:modec}
\scriptsize
\setlength{\tabcolsep}{3pt}
\begin{tabular*}{\textwidth}{@{\extracolsep{\fill}}l r r r r l l@{}}
\toprule
Instance & LB & UB & A/B & Mode C & Status & Method \\
\midrule
\multicolumn{7}{@{}l}{\emph{Job shop}} \\
\texttt{tai\_100\_1000\_1} & 522{,}298 & 533{,}080 & 525{,}515 & \textbf{525{,}389} & improved & compaction of the B schedule \\
\texttt{tai\_100\_1000\_4} & 527{,}101 & 534{,}746 & 528{,}156 & \textbf{527{,}101} & closed & corridor search, propagated windows \\
\texttt{tai\_100\_1000\_5} & 517{,}728 & 529{,}580 & 525{,}865 & \textbf{525{,}657} & improved & compaction of the B schedule \\
\texttt{tai\_100\_1000\_6} & 522{,}907 & 534{,}969 & 524{,}812 & \textbf{524{,}604} & improved & compaction of the B schedule \\
\texttt{tai\_100\_1000\_7} & 522{,}537 & 534{,}974 & 527{,}000 & \textbf{526{,}771} & improved & compaction of the B schedule \\
\texttt{tai\_100\_1000\_8} & 526{,}428 & 535{,}757 & 528{,}473 & \textbf{527{,}280} & improved & modified corridor search, compaction \\
\texttt{tai\_100\_1000\_9} & 528{,}097 & 536{,}993 & 529{,}470 & \textbf{528{,}097} & closed & corridor search, propagated windows \\
\texttt{tai\_100\_1000\_10} & 521{,}766 & 529{,}918 & 526{,}000 & \textbf{525{,}782} & improved & compaction of the B schedule \\
\midrule
\multicolumn{7}{@{}l}{\emph{Behnke--Geiger}} \\
\texttt{behnke14} & 221 & 230 & 232 & \textbf{228} & improved & new local search \\
\texttt{behnke32} & 220 & 224 & 223 & \textbf{220} & closed & CP-SAT at a fixed makespan \\
\texttt{behnke34} & 219 & 223 & 223 & \textbf{219} & closed & CP-SAT at a fixed makespan \\
\texttt{behnke36} & 381 & 388 & 389 & \textbf{386} & improved & new local search \\
\texttt{behnke38} & 386 & 389 & 391 & \textbf{387} & improved & CP-SAT at a fixed makespan \\
\texttt{behnke39} & 384 & 389 & 392 & \textbf{387} & improved & new local search \\
\texttt{behnke40} & 415 & 419 & 417 & \textbf{415} & closed & CP-SAT at a fixed makespan \\
\texttt{behnke52} & 210 & 212 & 214 & \textbf{211} & improved & CP-SAT at a fixed makespan \\
\texttt{behnke55} & 221 & 223 & 222 & \textbf{221} & closed & CP-SAT with a makespan cap \\
\texttt{behnke56} & 384 & 390 & 388 & \textbf{384} & closed & CP-SAT at a fixed makespan \\
\texttt{behnke57} & 385 & 390 & 389 & \textbf{386} & improved & new local search \\
\texttt{behnke58} & 392 & 397 & 398 & \textbf{394} & improved & new local search \\
\texttt{behnke60} & 397 & 402 & 400 & \textbf{397} & closed & CP-SAT at a fixed makespan \\
\bottomrule
\end{tabular*}
\end{table}

\end{document}